\documentclass[12pt,reqno]{amsart}

\usepackage{verbatim}

\usepackage{amsmath, amssymb}
\usepackage[latin1]{inputenc}

\newtheorem{theorem}{Theorem}[section]
\newtheorem{lemma}[theorem]{Lemma}
\newtheorem{prop}[theorem]{Proposition}
\newtheorem{cor}[theorem]{Corollary}

\newcounter{defn}

\newtheorem{remark}{Remark}[section]

\newtheorem{definition}[remark]{Definition}
\newtheorem{example}[remark]{Example}

\newcommand{\sect}{\vspace{3mm} \setcounter{equation}{0} \setcounter{remark}{0} \section}

\newcounter{num}

\newcommand{\slim}{\mbox{\rm s-}\lim}
\newcommand{\w}[1]{\langle {#1} \rangle}
\newcommand{\pf}{\noindent {\bf Proof. \hspace{2mm}}}
\newcommand{\ef}{ \hfill $ \blacksquare $ \vspace{3mm}}
\newcommand{\be}{\begin{equation}}
\newcommand{\ee}{\end{equation}}
\newcommand{\bea}{\begin{eqnarray}}
\newcommand{\eea}{\end{eqnarray}}
\newcommand{\beas}{\begin{eqnarray*}}
\newcommand{\eeas}{\end{eqnarray*}}
\newcommand{\ep}{{\epsilon}}
\newcommand{\f}{\frac}

\renewcommand{\d}{{\rm d}}

\newcommand{\supp}{{\rm supp}}

\renewcommand{\Re}{{\rm Re}\,}
\renewcommand{\Im}{{\rm Im}\,}

\newcommand{\mand}{\text{ and }}
\newcommand{\mfor}{\text{ for }}

\newcommand{\mor}{\text{ or }}

\newcommand{\bC}{{\mathbb C}}
\newcommand{\bR}{{\mathbb R}}
\newcommand{\bN}{{\mathbb N}}

\newcommand{\bS}{{\mathbb S}}

\newcommand{\vA}{{\mathcal A}}

\newcommand{\vD}{{\mathcal D}}

\newcommand{\vL}{{\mathcal L}}

\newcommand{\vR}{{\mathcal R}}

\newcommand{\vV}{{\mathcal V}}

\newcommand{\vS}{{\mathcal S}}

\renewcommand{\ker}{\rm Ker }

\newcommand{\ran}{\mbox{\rm Rank }}
\newcommand{\range}{\mbox{\rm Range }}

\def\p{\partial}

\def\f{\frac}

\def\la{\lambda}

\def\vp{\varphi}
\def\O{\Omega}

\begin{document}

\title[Gevrey estimates of the resolvent]
{ Gevrey estimates of  the resolvent and sub-exponential time-decay of solutions}
\author{ Xue Ping WANG}
\address{
Département de Mathématiques\\
UMR 6629 CNRS\\
Université de Nantes \\
44322 Nantes Cedex 3  France \\
E-mail: xue-ping.wang@univ-nantes.fr }
\date{\today}

\begin{abstract}
In this article, we study a class of non-selfadjoint Schr\"odinger operators $H$ which are perturbation of some model operator $H_0$  satisfying a weighted coercive assumption.  For the model operator $H_0$, we prove that the derivatives of the resolvent  satisfy some Gevrey estimates at threshold zero. As application, we establish large time expansions  of semigroups $e^{-tH}$ and $e^{-itH}$ for $t>0$ with subexponential time-decay estimates on the remainder, including possible presence of zero eigenvalue and  real resonances.
\end{abstract}

\maketitle

\tableofcontents

\sect{Introduction}

This work is concerned with time-decay of semigroups  $e^{-tH}$ and $e^{-itH}$ as $t \to +\infty$ where 
 $H = -\Delta + V(x)$ is a compactly supported perturbation of some model operator $H_0 = -\Delta + V_0(x)$ with a complex-valued potential  
 $V_0(x)= V_1(x) -i V_2(x)$ verifying a weighted coercive condition (see (\ref{ass2})). Large time behavior of solutions is closely related to low-energy spectral properties of operators.
  There are many works on low-energy spectral analysis and large time asymptotics  for selfadjoint Schr\"odinger operators $-\Delta + V(x)$ with a real-valued decreasing potential $V(x)$ verifying 
\be
|V(x)|\le C\w{x}^{-\rho}, \quad x\in \bR^n,
\ee
for some $\rho>0$, where $\w{x} = (1 + |x|^2)^{\f 1 2}$. Here we only mention \cite{bo,jk} for quickly decaying potentials ($\rho> 2)$, \cite{w1} for critically  decaying potentials ($\rho =2$) under an assumption  of Hardy inequality for the model operator 
and \cite{sw} in one-dimensional case when this Hardy condition is not satisfied. For slowly decreasing potentials ($0<\rho <2$), there are works of \cite{FS} when the potential is negative and \cite{n,y1,y2} when it is globally positive. When $\rho \ge 2$, threshold zero may be an eigenvalue and/or a resonance and for critically decaying potentials, threshold resonance may appear in any space dimension with arbitrary multiplicity. For slowly decreasing
potentials ($0<\rho <2)$, threshold resonance is absent and low-energy spectral analysis has not yet been done in presence of zero eigenvalue. For non-selfadjoint Schr\"odinger operators, we can mention works \cite{sai,ky} on the limiting absorption principle- and \cite{god,w2} on dispersive estimates. In \cite{god},  absence of real resonances is assumed and in \cite{w2} only dissipative operators are considered. In the later case, positive resonances may exist but outgoing positive resonances (see Definition \ref{def2.1}) are absent due to the dissipative condition on potential.
\\

Known results related to the topic studied in this work concern mainly  selfadjoint Schr\"odinger operators with globally positive and slowly decreasing potentials.  Let $H_0= -\Delta + V_0(x)$ with a real positive potential $V_0$ satisfying for some constants  $\mu \in ]0,1[$ and $c_1, c_2 >0$
 \be
 c_1\w{x}^{-2\mu} \le  V_0(x) \le c_2 \w{x}^{-2\mu}, \quad x\in \bR^n.
 \ee
 Under some additional conditions on $V_0$, it is known (\cite{n,y1,y2})) that the spectral measure $E_0'(\lambda)$ of $H_0$
 is  smooth at $\lambda =0$ and  satisfies for any $N \ge 0$ 
 \bea \label{1.5}
 \|E_0'(\lambda)\|_{L^2_{\rm comp} \to  L^2_{\rm loc}} &= & O(|\lambda|^N), \quad \lambda \to 0, \\
 \| e^{-tH_0}\|_{L^2_{\rm comp} \to  L^2_{\rm loc}} &=& O(e^{-ct^\beta})
 \eea
 where $\beta = \f{1-\mu}{1+\mu}$ and $c$ is some positive constant.
 In one dimensional case,  if $V_0(x)$ is in addition analytic, D. Yafaev (\cite{y1}) proves that
 \be \label{1.6}
  \|e^{-it H_0}\|_{L^2_{\rm comp} \to  L^2_{\rm loc}} = O(e^{-c|t|^\beta}), \quad |t|\to +\infty.
 \ee
  The proof given in \cite{y1} 
is based on explicit construction of solutions to Schr\"odinger equation in one dimensional case which is not available in higher dimensions. 
Another related topic is  return to equilibrium of Fokker-Planck operator with a positive, sublinearly increasing potential. In \cite{w3}, it is conjectured that the convergence rate in this situation should be subexponential in time. 
We learned from F. Bolley that there exists  probabilistic approaches to this problem. See the lecture notes of P. Cattiaux \cite{pc} for an overview. While polynomially decaying remainder estimate is proved in \cite{dfg}, the subexponential remainder estimate (\ref{eq1.6}) is proved in a recent work of T. Li and Z. Zhang (\cite{lz}).  Note that M. Klein and J. Rama (\cite{kr}) also studied Gevrey estimates in a different context to analyze large time evolution of quantum resonance states. \\

 In this article we are mainly interested in  non-selfadjoint Schr\"odinger operators $H=-\Delta +V(x)$, although  Gevrey estimates of the model resolvent at threshold are proved for a class of second order elliptic operators. Typically the potential in the model operator  $H_0= -\Delta + V_0(x)$ is of the form: $V_0(x) = V_1(x) -i V_2(x)$ with $V_1, V_2 \ge 0$  verifying
\be
c\w{x}^{-2\mu}\le  V_1(x) + V_2(x) \le C \w{x}^{-2\mu}, x\in \bR^n,
\ee
for some constants $c, C>0$ and $\mu \in ]0, 1[$. We assume throughout this paper that
\be 
\mbox{ $W= V-V_0$ is a bounded, compactly supported function.} 
\ee
If $V_0$ belongs to the class of  potentials $\vA$ (see Definition \ref{A}), one of the results proved  in this work for the Schr\"odinger semigroup  $e^{-itH}$  is an asymptotic expansion  of the form
\be \label{eq0.13}
 \| \chi  (e^{-itH} -\sum_{\lambda \in \sigma_d(H)\cap \overline{\bC}_+} e^{- itH} \Pi_{\lambda} -\Pi_0(t) - \sum_{\nu \in r_+(H)} e^{-it \nu} P_\nu(t) )\chi\|_{\vL(L^2)} \le C_\chi  e^{-ct^{\f{1-\mu}{1+\mu}}}  \quad t>0. 
\ee
Here $\chi \in C_0^\infty(\bR^n)$, $c>0$ is independent of $\chi$,  $r_+(H)$ is the set of outgoing positive resonances of $H$ (see  Definition \ref{def2.1}), $\Pi_\lambda$ is the Riesz projection of $H$ associated with $\lambda$,
 $\Pi_0(t)$ and $P_\nu(t)$   are operators of finite rank depending  polynomially on $t$.  See Theorem \ref{th1.4} for more precision on conditions and results. Note that if $\Im V \le 0$, then 
$r_+(H) =\emptyset$ and that there exists $V\in C_0^\infty(\bR^n)$ with $\Im V \ge 0$ such that  $r_+(-\Delta + V) \neq \emptyset$ (see \cite{w2} for an example of incoming positive resonance with $V\in C_0^\infty$ and $\Im V \le 0$).
\\

To prove (\ref{eq0.13}), we use analytic deformation of $H$ outside some sufficiently large ball in $\bR^n$ and prove the existence of a curve $\Gamma $ 
in the lower half complex plane, intersecting the real axis only at point $0$,  such that  the above this curve, the meromorphic extension of cut-off resolvent
$\chi R(z)\chi$ from $\bC_+$  has  at most a finite number of poles and those located in $]0, +\infty[$ are precisely outgoing positive resonances. In particular, zero is not an accumulation point of quantum resonances of $H$ located in that region. 
Under some conditions, we compute the resolvent expansion at threshold  in presence of zero eigenvalue   and prove the Gevrey estimates  for the remainder. Then
(\ref{eq0.13}) is deduced by representing $\chi e^{-itH} \chi $ as  sum of some residue terms and a Cauchy integral of the cut-off resolvent on $\Gamma$. The subexponential time-decay estimates is obtained from the Gevrey estimates on the remainder by expanding it at threshold up to some order $N$ appropriately chosen according to $t>0$ .
\\

Real resonances, called spectral singularity  by J. Schwartz in \cite{sch} in more general framework, are the main obstacle to study spectral properties of non-selfadjoint Schr\"odinger operators near positive real half-axis.  Up to now, one only knows that  real resonances form  bounded set with Lebesgue measure zero (\cite{sai,sch}). To study spectral properties  of non-selfadjoint Schr\"odinger operators near positive real half-axis, one usually  either supposes the absence of real resonances or uses some kind of exponential-type decay on potentials. In this work we prove for some classes of potentials with analyticity condition, real resonances are at most a countable set with
zero as the only possible accumulation point. If in addition the weighted coercive condition is satisfied, then outgoing positive resonances are at most a finite set. Each outgoing positive resonance is a pole of some meromorphic extension of the resolvent from the upper half complex plan, hence contributes to large time asymptotics of solutions as $t\to +\infty$.
\\
 
A technical task in this work is Gevrey estimates on various remainders at threshold zero. To establish  Gevrey estimates on the resolvent of model operator $H_0$ at threshold, we prove firstly an energy estimate uniformly on some parameter $s\in \bR$. This kind of estimate for fixed $s$ is already proved by D. Yafaev in \cite{y2}. The uniformity on $s\in \bR$ is crucial for us, because it allows to control  norms of the resolvent in weighted spaces with respect to several parameters (Theorem \ref{th3.2}), from which we deduce   Gevrey estimates on the model resolvent at threshold (Theorem \ref{th1.1}). To estimate remainders  in the asymptotic expansions of $(H-z)^{-1}$ near $0$, we make use of Theorem \ref{th1.1} for the model operator and  operations of operator-valued functions in Gevrey classes. \\

The organisation of this paper is as follows.  In Section 2, we introduce  conditions and definitions used and  state  main results obtained in this work.
Sections 3 and 4 are devoted to the analysis of the model operator $H_0 = -\Delta + V_0(x)$ verifying the weighted coercive condition (\ref{ass2}). In Section 3,  we prove Gevrey estimates of the model resolvent at threshold (Theorem \ref{th1.1}). We firstly establish a uniform  energy estimate which allows to control the growth of powers of the resolvent at threshold in weighted spaces. Then Theorem \ref{th1.1} is deduced by an appropriate induction. In Section 4, we begin with evaluating the numerical range of $H_0$ and prove resolvent estimates for $H_0$ on the left of a curve located in the right half complex-plane. An estimate like (\ref{eq1.8}) is proved for $H_0$.
We prove the absence of complex eigenvalues in some domaine near zero for a class of  model operators $H_0$. The subexponential time-decay of $e^{-itH_0}$ is studied in subsection 4.3 when potential $V_0$ belongs to some analytic class $\vA$ introduced in Section 2. We show that there exists a contour located in the lower half complex plane passing by $0$ on which the cut-off resolvent $\chi(H_0-z)^{-1}\chi$ is uniformly bounded and that  there are no quantum resonances and eigenvalues of  $H_0$ in a sector below the positive half real axis. 
(\ref{eq1.9}) for $H_0$ is then obtained by deforming the integral contour into the lower half complex plane.  
Compactly supported perturbations of the model operator $H_0$ are studied in Sections 5 and 6.  In Section 5, we  analyze properties of real resonances, study the low-energy resolvent expansion for $H$ and prove Theorems \ref{th1.2} and \ref{th1.3}. Since the method of low-energy spectral analysis used in the proof of Theorem \ref{th1.3} is well known for selfadjoint operators, we only emphasize upon Gevrey estimates on remainders. Subexponential time-decay estimates in Theorems \ref{th1.2} and \ref{th1.3} are derived from low-energy resolvent expansion by the same method as that used for $H_0$.
 Finally in Section 6, we study the case of threshold eigenvalue for non-selfadjoint Schr\"odinger operators and prove Theorem \ref{th1.4}.
 In order to obtain more precisions when zero eigenvalue is geometrically simple, We firstly establish a representation formula for the Riesz projection $\pi_1$  associated to the compact operator $G_0W$ with eigenvalue $-1$  and then use Grushin method to compute the leading term of the resolvent.  The Gevrey estimates on remainders can be proved as in Section 5 and hence the details are omitted in Section 6.
Results related to  the model operator $H_0$ are announced in \cite{w4}.
\\
 
 {\noindent \bf Notation}.  We denote $H^{r,s}$, $r \ge 0, s\in\bR$ the weighted Sobolev space of order $r$  with  the weight $\w{x}^{s}$ on $\bR^n$:
\[
H^{r,s} =\{u\in \vS'(\bR^n); \|u\|_{r,s} =\|\w{x}^{s}(1-\Delta)^{\f r 2}u \|_{L^2} <\infty\}.
\]
 For $r <0$, $H^{r, s}$ is defined as dual space of $H^{-r,-s}$ with dual product identified with the scalar product 
$\w{\cdot, \cdot}$ of $L^2(\bR^n)$.   Denote $H^{0,s} = L^{2,s}$. $\vL({r,s}; {r',s'})$ stands for the space of continuous linear operators from $H^{r,s}$ to $H^{r',s'}$.  If $(r,s) =(r',s')$, we denote $\vL(r,s) =\vL({r,s}; {r',s'})$. Unless mentioned explicitly, $\| \cdot\|$ denotes norm in $L^2(\bR^n)$ or in $\vL(L^2)$ when no confusion is possible. $\bC_\pm$ denote respectively the upper and the lower open half-complex plane  and $\overline{C}_\pm$ their closure. \\

\sect{Statement of results}

We shall prove Gevrey estimates of the resolvent at threshold for a class  second order elliptic operators satisfying a weighted coercive condition. Let
 \be \label{e1.1} 
 H_0 = -\sum_{i,j =1}^n \partial_{x_i} a^{ij}(x) \partial_{x_j} +  \sum_{j =1}^n  b_j(x) \partial_{x_j}  + V_0(x), 
 \ee 
 where $a^{ij}(x)$, $b_j(x)$ and $V_0(x)$  are  complex-valued measurable functions. Suppose that
  $a^{ij}, b_j \in C_b^1(\bR^n)$  and that there exists $c>0$
 such that 
 \be \label{e1.2}
 \Re (a^{ij}(x)) \ge c I_n, \quad \forall x \in \bR^n.
 \ee
Assume that $V_0$ is relatively bounded with respect to $-\Delta$ with relative bound zero,  $\Re H_0 \ge 0$ and that there exists some constants $0<\mu<1$ and $c_0>0$ such that
\bea \label{ass2} 
& &|\w{H_0 u, u}| \ge  c_0 ( \|\nabla u\|^2 + \|\w{x}^{-\mu}u\|^2), \quad \mbox{for all } u \in H^2,  \\[3mm]
& & \sup_x |\w{x}^\mu b_j (x) | < \infty, \quad j = 1, \cdots, n. \label{ass2b}
\eea 
Condition (\ref{ass2}) is called weighted coercive condition.\\

\begin{remark} If $H_0 = -\Delta + V_0(x)$ with $V_0(x) = V_1(x)-i V_2(x)$ with $V_j $ real. Assume that
$-\alpha \Delta + V_1(x) \ge v_1(x) \ge 0$  for some  $\alpha \in ]0, 1[$. If $V_2\ge 0$ is such that for some $c>0$
\be \label{v1}
 v_1(x) + V_2(x) \ge c\w{x}^{-2\mu}, \quad x\in \bR^n.
\ee 
then  the weighted coercive condition (\ref{ass2}) is satisfied. If $V_1(x)$ is globally positive and slowly decaying ({\it i. e.}  $V_1(x) \ge c\w{x}^{-2\mu}$ for some $\mu \in ]0,1[$ and $c>0$), then 
(\ref{ass2}) is satisfied by $H_0 = -\Delta + V_1(x) - iV_2(x)$ for any real function $V_2$ which is $-\Delta$-bounded with relative bound zero.
\end{remark}

Note that when we study Schr\"odinger operators $H_0= -\Delta + V_0$ by technics of analytic dilation or  analytic deformation,  if $H_0$ verifies (\ref{ass2}), the analytically dilated or distorted operators  obtained from $H_0$ are of the form (\ref{e1.1}) and still satisfy (\ref{ass2}) if the dilation or distortion parameter is small.  The condition that $V_0(x)$ is $-\Delta$-bounded 
with relative bound zero allows to include a class of $N$-body potentials.\\

Under the assumptions \ref{e1.2}, \ref{ass2} and \ref{ass2b},  one can show that $H_0$ is bijective from $D(H_0)= H^2(\bR^n)$ to $R(H_0)$ and $R(H_0)$ is dense in $L^2(\bR^n)$. Let $G_0 : R(H_0) \to D(H_0)$ be the algebraic inverse of $H_0$. Denote $L^{2,s}= L^2(\bR_x^n; \w{x}^{2s} dx)$ and $\vD= \cap_{s\in \bR} L^{2,s}$.
 Then  $G_0(\vD) \subset \vD$ and $G_0$ is a densely defined, continuous from $R(H_0)\cap L^{2,s}$ to $L^{2, s-2\mu}$ for any $s \in \bR$.
 (See Lemma \ref{lem3.1}).  To simplify notation, we still denote by $G_0$ its continuous extension by density so that $G_0$ is regarded as a bounded operator from $L^{2,s}$ to $L^{2, s-2\mu}$. Consequently for any 
$N\in \bN$,  $G_0^N : L^{2, s} \to L^{2, s-2\mu N}$ is well defined for any $s\in \bR$. 
Let $R_0(z) = (H_0-z)^{-1}$ for $z\not\in \sigma(H)$. Since $\Re H_0 \ge 0$ on $L^2$, one has
\[
\slim_{z\in \Omega(\delta), z\to 0} R_0(z) = G_0 
\]
as operators from $L^{2,s}$ to $L^{2, s-2\mu}$, where $\Omega (\delta)=\{z;  \f {\pi} 2 + \delta < \arg z  < \f {3\pi} 2 - \delta\}$ for some $\delta>0$.\\

\begin{theorem}\label{th1.1} Assume the conditions  (\ref{e1.1})-(\ref{ass2b}).
The following estimates hold. \\

(a).  For any $a>0$,  there exists $C_a>0$ such that
\be \label{eq1.6}
 \|e^{-a\w{x}^{1-\mu}}G_0^N \|+  \|G_0^N e^{-a\w{x}^{1-\mu}}\|\le C_a^N N^{\gamma N}, \forall N.
\ee

(b).  There exists some constant $ C>0$ such that $\forall \chi\in C_0^\infty(\bR^n)$, one has for some $C_\chi >0$
 \be \label{eq1.7}
 \|\chi(x)G_0^N \| + \|G_0^N \chi(x)\|\le C_\chi C^N N^{\gamma N}, \forall N.
\ee
Here $\gamma = \f{2\mu}{1-\mu}$.
\end{theorem}

 Since one has at least formally
\[
\f{d^N}{dz^N} R_0(z)|_{z=0} = N! G_0^{N+1},
\]
 Theorem \ref{th1.1} says that derivatives of the resolvent of $H_0$ at threshold satisfies the Gevrey estimates of order $\sigma = 1+\gamma$.
To study large time behavior  of semigroups, we introduce two classes of potentials $\vV$ and $\vA$.

\begin{definition} \label{V}
Denote $\vV$  the class of  complex-valued potentials   $V_0$ such that
\be
\label{ass2c}
\mbox{ $V_0$ is $-\Delta$-compact and (\ref{ass2}) is satisfied for some $\mu \in ]0, 1[$. } 
\ee
and
\be
\label{vV}
 \Re H_0  \ge -\alpha \Delta \mbox{ and } |\Im V_0(x)| \le C \w{x}^{-2\mu'}
\ee
for some constants $\alpha, \mu', C>0$. 
\end{definition}

Results for the heat semigroup $e^{-tH}$ will be proved for model potentials $V_0\in \vV$.
To study the Schr\"odinger semigroup $e^{-itH}$ we will use both technics of analytical dilation and analytical deformation, hence need the analyticity of potentials. 

\begin{definition} \label{A}
Let $\vA$ denote the class of complex-valued  potentials
 $V_0(x) = V_1(x)- iV_2(x) $ for $x\in \bR^n$ with $n \ge 2$  such that $V_0\in \vV$  satisfies the estimate (\ref{ass2}) for some  
 $\mu\in ]0,1[$ verifying
\be
\mbox { $0<\mu <\f 3 4$ if $n=2$ and $0<\mu < 1$
 if $n \ge 3$;} 
 \ee
   and that $V_1$ and $V_2$ are  dilation analytic (\cite{ac})  and extend holomorphically into a complex region of the form 
 \[
 \Omega =\{ x\in \bC^n; |\Im x |< c |\Re x |\} \cup\{x\in \bC^n; |x| >c^{-1}\}
 \]
  for some $c >0$ and satisfy for some $c_j >0$ and $R\in[0, +\infty]$
 \bea
|V_j(x)| & \le & c_1 \w{\Re x}^{-2\mu}, x\in \Omega, \quad j =1,2, \label{e1.25} \\
  V_2(x) & \ge & 0 , \quad \forall x\in \bR^n, \label{e1.26}  \\
 x\cdot\nabla V_1 (x) &\le & -c_3 \f{x^2}{\w{x}^{2\mu +2}}, \quad x\in \bR^n \mbox{ with } |x| \ge R, \mbox{ and } \label{e1.27}
 \\
  V_2(x) &\ge & c_5 \w{x}^{-2\mu},
 \quad x\in \bR^n \mbox{ with } |x| < R.  \label{e1.28}
\eea 
\end{definition}

Remark that when $R =0$, (\ref{e1.27}) is a global virial condition on $V_1$ and (\ref{e1.28})
 is void; while if $R= +\infty$, no virial condition is needed on $V_1$, but (\ref{e1.28})
 is required on the whole space which means that the dissipation is strong. Potentials of the form
 \be
 V_0(x) = \f{c}{\w{x}^{2\mu}} - iV_2(x)
 \ee
 satisfy conditions (\ref{e1.25}-(\ref{e1.28}) with $R=0$, if  $V_2 \ge 0$ and $V_2$ is  holomorphic in $\O$ satisfying 
 $|V_2(x)| \le C \w{x}^{-2\mu}$ for $x\in \O$. Conditions (\ref{e1.27}) and  (\ref{e1.28}) are used to prove the non-accumulation of quantum resonances towards zero in some sector.\\

 For  $V_0 \in \vA$, one can study quantum resonances
 of  $H_0= -\Delta + V_0(x)$ by both analytical dilation or analytical deformation outside some compact (\cite{ac,hun,si}).
 We shall show that under the conditions (\ref{e1.26}), (\ref{e1.27}) and (\ref{e1.28}), there are no  quantum resonances
 of $H_0$ in a sector  below the positive real half-axis in complex plane. \\

Let  $V_0 \in \vV$ and $H_0 = -\Delta + V_0(x)$. Let $H= H_0 +W(x)$  be a compactly supported perturbation of $H_0$:  
$W \in L^\infty_{\rm comp} =\{u \in L^\infty(\bR^n),\supp \; u  \mbox{ compact }\}$.  It will be proved in Section 6 that $H$ has only at most a finite number of  discrete eigenvalues located on the left of a curve $\Gamma$ of the forme 
\[
\Gamma =\{z; \Re z \ge 0,  | \Im z|= C (\Re z)^{\mu'}\}
\]
for some constants $C, \mu' >0$ and that there exists a nice bound for the resolvent of $H_0$ on $\Gamma$.  Note that  zero may be an embedded eigenvalue, but it is never a resonance of $H$,  {\it i. e.},  if 
$u \in L^{2}(\bR^n; \w{x}^{2s} dx)\cap H^{2}_{{\rm loc}}(\bR^n)$ for some $s\in \bR$ such that $Hu=0$, then one can show that 
$u \in H^2(\bR^n)$. Complex eigenvalues of $H$ may accumulate to zero from the right side of $\Gamma$. 
Let $\sigma_d(H)$ ($\sigma_p(H)$, resp.)  denote the set of discrete eigenvalues of $H$ (the set of eigenvalues of $H$, resp.). 
\\

More subtle is the role of real resonances. Recall that if $V$ is of short-range, $\lambda >0$ is called resonance of $H=-\Delta + V(x)$
if the equation $Hu = \lambda u$ admits a non-trivial solution $u\in H^2_{\rm loc}(\bR^n)$ satisfying one of Sommerfeld radiation conditions:
\be \label{rad}
u(x) = \f{e^{\pm i\sqrt{\lambda} |x|}}{|x|^{\f{n-1} 2}} (a_\pm(\omega) + o(1)), \quad |x| \to \infty,
\ee
for some $a_\pm \in L^2(\bS^{n-1}),  a_\pm \neq 0$. $\lambda$ is called an outgoing (resp., incoming) positive resonance of $H$ if $u$ verifies (\ref{rad}) with sign $+$ (resp. with sign $-$). It is known that  if $V$ is real, then positives resonances are absent (\cite{ag}) and if 
$\Im V \le 0$, outgoing resonances are absent (\cite{roy0}). In this paper, we use a slight different definition for outgoing resonances. Let $U_0$ be a complex valued function such that
$(x\cdot \nabla_x)^jU_0$, $j=0,1,2$, are $-\Delta$-compact and $\Im U_0 \le 0$. Then for any $\la >0$ the boundary value of the resolvent
\be
(-\Delta + U_0 - (-\la + i0))^{-1} = \lim_{z\to \la, \Im z >0} (-\Delta + U_0 - (z)^{-1} 
\ee 
exists in $\vL(0, s; 0, -s)$ for any $s > \f 1 2$  and is H\"older-continuous for $\la >0$. See \cite{roy}. 

\begin{definition} \label{def2.1}  Let  $U(x)$ be a Lebesgue measurable function such  that $U(x) -U_0(x)$ is bounded and of short-range on $\bR^n$.
$\lambda >0$ is called outgoing resonance of $-\Delta + U$ if $-1$ is an eigenvalue of the compact operator 
$(-\Delta + U_0 - (-\la + i0))^{-1}(U-U_0)$ in $L^{2,-s}$ for $s>\f 1 2$ and sufficiently close to $\f 1 2$.
Denote $r_+(-\Delta + U)$ the set of outgoing  resonances of $H$. For $\la \in r_+(H)$, define $m_+(\la)$ as the algebraic multiplicity of eigenvalue $-1$ of $(-\Delta + U_0 - (-\la + i0))^{-1}(U-U_0)$. 
Similarly if $\Im U_0 \ge 0$, one can define the set of incoming positive resonances $r_-(-\Delta+ U)$  and $m_-(\la)$ for $\la \in r_-(-\Delta+ U)$.
\end{definition}

If $U$ is of short-range, then our definition coincides with the usual one. 
 In the case that zero is not an eigenvalue of non-selfadjoint Schr\"odinger operator $H=H_0 + W(x)$,  we prove the following \\

\begin{theorem}\label{th1.2}
Assume that zero is not an eigenvalue of $H$.  \\

(a).  Let $V_0 \in \vV$. For any $a>0$ there exist $c_a, C_a >0$ such that 
\be \label{eq1.8}
 \| e^{-a\w{x}^{1-\mu}} (e^{-tH} -\sum_{\lambda \in \sigma_d(H), \Re \lambda \le 0} e^{-tH} \Pi_{\lambda})\| \le C_a  e^{-c_a t^{\beta}}  \quad t>0, 
\ee
where 
\be \label{beta}
\beta =  \f{1-\mu}{1+\mu}.
\ee

(b).  Let  $V_0\in \vA$.   Then the set of outgoing resonances $r_+(H)$ of $H$  is at most finite. There exists some constant $c>0$ such that 
 for any $\chi \in C_0^\infty(\bR^n)$ one has
\be \label{eq1.9}
 \| \chi  (e^{-itH} -\sum_{\lambda \in \sigma_d(H)\cap \overline{\bC}_+ }e^{- itH} \Pi_{\lambda} - \sum_{\nu \in r_+(H)} e^{-it \nu} P_\nu(t))\chi\| \le C_\chi  e^{-c\;t^{\beta}}  \quad t>0, 
\ee
Here $\Pi_\lambda$ denotes the Riesz projection associated with the discrete eigenvalue $\lambda$ of $H$ and $P_\nu(t)$ is an operator depending polynomially on  $t$  with coefficients of rank not exceeding $m_+(\nu)$. 
\end{theorem}

Consider now the case that zero is an eigenvalue of $H$. If $H$ is selfadjoint, $H$ has only a finite number of negative  eigenvalues and both positive eigenvalues and positive resonances are absent.  We can apply the known  method in threshold spectral analysis for selfadjoint operators to compute low-energy expansion of the resolvent. Theorem \ref{th1.1} allows to estimate remainders in Gevrey spaces and to prove the following\\

\begin{theorem}\label{th1.3}  Assume that zero is an eigenvalue of $H$ and that both $H$ and $H_0$ are selfadjoint.
\\

(a). If $V_0\in \vV$, then for any $a>0$,  there exist some constants $c_a, C_a >0$ such that 
\be \label{eq1.10}
 \| e^{-a\w{x}^{1-\mu}} (e^{-tH} -\sum_{\lambda \in \sigma_p(H)} e^{-t\lambda} \Pi_{\lambda})\| \le C_a  e^{-c_a t^{\beta}}  \quad t>0, 
\ee

(b).  Let  $V_0\in \vA$.  Then there exists some constant $c>0$ such that 
 for any $\chi \in C_0^\infty(\bR^n)$, one has
\be \label{eq1.11}
 \| \chi  (e^{-itH} -\sum_{\lambda \in \sigma_p(H)} e^{- it\lambda} \Pi_{\lambda})\chi\| \le C_\chi  e^{-c\; t^{\beta}}  \quad t>0, 
\ee
Here $\Pi_\lambda$ denotes the orthogonal eigenprojection of $H$  associated with eigenvalue $\lambda$ of $H$.
\end{theorem}

Theorem \ref{th1.3} can be applied to a class of Witten Laplacians for which  zero is an eigenvalue embedded in the continuous spectrum which is equal to $[0, +\infty[$. Our result is new  concerning the Schr\"odinger semigroup $e^{-itH}$ in presence of zero eigenvalue. For the heat semigroup,  there are results obtained by method of Markov processes (\cite{pc, dfg}). \\

The case of zero eigenvalue  in non-seladjoint case is more difficult. There does not yet exit general method  to treat this kind of problem.    In our case $G_0W$ is compact on $L^2(\bR^n)$ and  zero is an eigenvalue of $H$ if and only if $-1$ is an eigenvalue of $G_0W$ and zero eigenvalue of $H$, if it does exist, is of finite geometrical multiplicity. 
Let $m$ denote the algebraic multiplicity of eigenvalue  $-1$  of $G_0W$. 
Then one can show that there exists some numerical Gevrey function $\omega(z) $  such that for $z$ near zero, $z\in \sigma_d(H)$ if and only if $\omega(z) =0$. (See Proposition \ref{prop7.1}). In addition, $\omega(z)$ admits an asymptotic expansion of any order in powers of $z$: there exist some constants $\omega_j \in \bC$, $j \in \bN$,  such that 
\be
\omega(z)  = \omega_1 z + \cdots  + \omega_N z^N + O(|z|^{N+1}), 
\ee
for $z$ near $0$ and $\Re z <0$ and for any $N \in \bN^*$.

\begin{theorem}\label{th1.4} 1. Assume that zero is an  eigenvalue of $H$ and that there exists  some $\omega_k \neq 0$ such that
\be \label{assk}
\omega(z) = \omega_k z^k +  O(|z|^{k+1}), 
\ee
for $z$ near $0$ and $\Re z <0$. Then  the following results hold.\\

1a. If $V_0\in \vV$, then for any $a>0$,  there exists some constants $c_a, C_a >0$ such that 
\be \label{eq1.12}
 \| e^{-a\w{x}^{1-\mu}} (e^{-tH} -\sum_{\lambda \in \sigma_d(H), \Re \lambda \le 0} e^{-t H }\Pi_{\lambda} -\Pi_0(t))\| \le C_a  e^{-c_a t^{\beta}}  \quad t>0, 
\ee

1b.  Let  $V_0\in \vA$.   Then the set of outgoing resonances $r_+(H)$ of $H$ is at most finite and there exist $c>0$ such that
 for any $\chi \in C_0^\infty(\bR^n)$, 
\be \label{eq1.13}
 \| \chi  (e^{-itH} -\sum_{\lambda \in \sigma_d(H)\cap \overline{\bC}_+} e^{- itH} \Pi_{\lambda} -\Pi_0(t) - \sum_{\nu \in r_+(H)} e^{-it \nu} P_\nu(t) )\chi\| \le C_\chi  e^{-c\; t^{\beta}}  \quad t>0, 
\ee
Here  $\Pi_\lambda$ and $P_\nu(t)$ have  the same meaning as in Theorem \ref{th1.2} (b) and $\Pi_0(t)$ is polynomial in $t$ of the form
\be
\Pi_0(t) = \sum_{j=0}^{k-1} t^j\Pi_{0,j}
\ee
where $\Pi_{0, j}$, $ 0 \le j \le k-1$ is an operator of rank not exceeding $m$, $m$ being the algebraic multiplicity of $-1$ as eigenvalue of $G_0W$. \\

2. Assume that  zero eigenvalue of $H$ is geometrically simple. \\

2a.  If (\ref{assk}) is satisfied, then $ \Pi_{0,k-1} $ is of rank one,  given by
 \be
 \Pi_{0,k-1} =\w{\cdot, J\psi_0}\psi_0
 \ee
  for some eigenfunction $\psi_0$ associated with zero eigenvalue of $H$. Here $J$ is the complex conjugaison 
  $J : f(x) \to \overline{f(x)}$.
  \\
  
2b. If there exists an associated eigenfunction $\varphi_0$ such that
\be \label{1.41}
\int_{\bR^n} (\varphi_0(x))^2 dx =1,
\ee
then Condition (\ref{assk}) is satisfied with $k=1$ and one has
\be
\Pi_0(t) = \Pi_{0,0} =\w{\cdot, J\vp_0}\vp_0.
\ee
\end{theorem}

Condition (\ref{assk}) is similar to that used in \cite{sch} to study spectral properties of some non-selfadjoint operators in presence of spectral singularities. This kind of conditions can be checked at positive resonances under some analyticity assumptions (see Remark \ref{rmk7.1}. But it is not clear if it can be satisfied at threshold even if the potential is analytic.  The assumption that zero is a geometrically simple eigenvalue of $H$ implies that there is only one Jordan block of the compact operator $G_0W$ associated with eigenvalue $-1$. This allows to construct explicitly a representation of the associated Riesz projection and to compute the leading term. The method developed in the proof of Theorem \ref{th1.4} is general and applies to some other situations. See Remarks \ref{rmk7.1} and \ref{rmk7.2}. In Section 6, we give an example such that (\ref{1.41}) is satisfied. \\

\sect{Gevrey estimates of the model resolvent at threshold}

 The starting point of our Gevrey estimates of the resolvent of $H_0$ is a uniform {\it a priori} energy estimate for the model operator $H_0$. In the sequel, we need to apply
 this kind of energy estimates to the Schr\"odinger operator $-\Delta + V_1(x)-i V_2(x)$ and
 to its analytically dilated or distorted versions as well.
For this purpose, we  begin with a slightly more general setting where  $H_0$ is a second order elliptic differential operator of the form
 \be \label{e2.1} 
 H_0= -\sum_{i,j =1}^n \partial_{x_i} a^{ij}(x) \partial_{x_j} +  \sum_{j =1}^n  b_j(x) \partial_{x_j}  + V_0(x)
 \ee 
 satisfying conditions (\ref{e1.2}), (\ref{ass2}) and (\ref{ass2b}). 
\\

Denote $b =(b_1, \cdots, b_n)$ and
\be \label{a}
 |a|_\infty =\max_{1\le i,j\le n}\sup_{x\in \bR^n} |a^{ij}(x)|, 
 \quad |b|_{\mu, \infty} =\max_{1\le j\le n}\sup_{x\in  \bR^n} |\w{x}^\mu b_{j}(x)|. 
 \ee 
 For $s \in \bR$, denote \be \varphi_s(x) = (1  + \f{|x|^2}{R_s^2})^{ s }, 
 \ee
 where $R_s =  M \w{s}^{\f 1{1-\mu}}$ with $M= M(c_0, |a|_\infty,  |b|_\infty)>1$ large enough, but
  independent of $s \in \bR$.  The uniformity in $s \in \bR$ in the following lemma is important for 
  Gevrey estimates  of the  model resolvent at threshold. \\

\subsection{A uniform  energy estimate}

\begin{lemma} \label{lem2.2} Let $H_0$ be given by (\ref{e2.1}).  Under the conditions (\ref{e1.2}), (\ref{ass2}) and (\ref{ass2b}) with $0<\mu <1$, there exist some constants $C, M >0$ depending only on $|a|_\infty$,  $|b|_{\mu,\infty}$ and $c_0$ given in (\ref{ass2}) such that 
\be \label{e2.11}
  \| \w{x}^{-\mu} \varphi_s(x) u\| +   \| \nabla (\varphi_s(x) u)\| \le  C \| \w{x}^{\mu} \varphi_s(x) H_0u \|
\ee 
for any $s\in\bR$ and $u  \in  H^2_{\rm loc}(\bR^n)$ with $ \w{x}^{s +\mu} H_0 u \in L^2$.
\end{lemma}
\pf We calculate $\w{u, \varphi_s^2 H_0 u}$ for $u \in C_0^\infty$: 
\bea 
\lefteqn{\w{u, \varphi_s^2 H_0 u}} \nonumber \\
 & = & \w{\varphi_s u, H_0(\varphi_s u)} +
 \w{\varphi_su, [\sum_{i,j =1}^n \partial_{x_i} a^{ij} \partial_{x_j}, \varphi_s] u}
  -\w{\varphi_s u,(b\cdot\nabla\varphi_s) u)}  \nonumber   \\
  &=& I + II + III,  \label{e2.12} 
\eea
where
\beas
I & = & \w{\varphi_s u, H_0(\varphi_s u)} \\
II &= & \w{\varphi_su, \sum_{i,j =1}^n \left((\partial_{x_i}\varphi_s) a^{ij} \partial_{x_j}u +
\partial_{x_i}( a^{ij} (\partial_{x_j}\varphi_s)u)\right)}  \\
III &=&  -\w{\varphi_s u,(b\cdot\nabla\varphi_s) u)}. 
\eeas 
Since $ \varphi_s  \partial_{x_j}u =   \partial_{x_j}(\varphi_s u) - (\partial_{x_j}\varphi_s) u$, one has \beas \lefteqn{|\w{\varphi_su, (\partial_{x_i}\varphi_s) a^{ij} \partial_{x_j}u +
\partial_{x_i}( a^{ij} (\partial_{x_j}\varphi_s))u)}|} \\
& = & |\w{(\partial_{x_i}\varphi_s)u,  a^{ij}  (\partial_{x_j}(\varphi_s u)
 -  (\partial_{x_j}\varphi_s) u)} + \w{\varphi_su,
\partial_{x_i}( a^{ij} (\partial_{x_j}\varphi_s)u)}|\\
& = & |\w{(\partial_{x_i}\varphi_s)u,  a^{ij}  (\partial_{x_j}(\varphi_s u)
-  (\partial_{x_j}\varphi_s) u)} - \w{\partial_{x_i}(\varphi_su), a^{ij} (\partial_{x_j}\varphi_s)u}| \\
&\le & |a|_\infty (\|(\partial_{x_i}\varphi_s)u\| (\|\partial_{x_j}(\varphi_s u)\|
 +  \|(\partial_{x_j}\varphi_s) u\|) + \|\partial_{x_i}(\varphi_su)\| \| (\partial_{x_j}\varphi_s)u\|)
\eeas The term $II$ in (\ref{e2.12}) can be bounded by \beas | II | &\le & |a|_\infty (\sum_{i=1}^n\|(\partial_{x_i}\varphi_s)u\| )
(\sum_{j=1}^n (2\|\partial_{x_j}(\varphi_s u)\| +  \|(\partial_{x_j}\varphi_s) u\|) ) \\
&\le &  n^2|a|_\infty \|(\nabla\varphi_s)u\| )(2 \|\nabla(\varphi_s u)\| +  \|(\nabla\varphi_s) u\|) ) \\
&\le &  n^2|a|_\infty  (\ep \|\nabla(\varphi_s u)\|^2 +  (1+ \f 1 \ep) \|(\nabla\varphi_s) u\|^2)
 \eeas 
 for any $\ep >0$. Clearly, $III$ verifies
 \be
 |III| \le |b|_{\mu, \infty} \|\w{x}^{-\mu}\varphi_s u\| \| (\nabla\varphi_s) u\|
 \le |b|_{\mu, \infty} (\ep   \|\w{x}^{-\mu}\varphi_s u\|^2
 +  \f 1 {4 \ep}\|(\nabla\varphi_s) u\|^2) 
 \ee
 Taking $\ep = \ep (c_0, |a|_\infty, |b|_{\mu, \infty})>0$ appropriately small where $c_0>0$ is given by (\ref{ass2}), 
 it follows from (\ref{ass2}) that
\bea \label{lowerb1}
|\w{u, \varphi_s^2 H u} | &\ge &|I| -| II| -|III|\\
&  \ge & \f{c_0}{2} (\|\nabla(\varphi_s u)\|^2 + \| \w{x}^{-\mu} \varphi_s(x) u\|^2) -  \w{u, W_s  u} \nonumber 
\eea 
where $W_s(x) = c_1 |\nabla \varphi_s|^2$ with $c_1>0$ some constant depending only
 on $c_0$, $|a|_\infty$ and $|b|_{\mu, \infty}$. One can check that
\beas
|\nabla \varphi_s|^2 & = & \f{4 s^2 x^2}{R_s^4 (1 + \f{x^2}{R_s^2})^2} (1 + \f{x^2}{R_s^2})^{2s} \\
& \le & \f{4s^2 x^2}{(R_s^2 +x^2)^2}  \varphi_s^2 \le \f{4s^2}{R_s^2  + x^2}  \varphi_s^2 
\eeas 
Since $R_s^2 + x^2 \ge 2^{-2\mu}R_s^{2(1-\mu)}\w{x}^{2\mu}$ and $R_s = M \w{s}^{\f 1 {1-\mu}}$,  $W(x)$ is bounded by 
\be 
0  \le W_s(x) \le \f{4 c_1 \w{s}^2}{R_s^2 + x^2}  \varphi_s^2\le \f{ 2^{2\mu} 4 c_1}{M^{2(1-\mu)} \w{x}^{2\mu}} \varphi_s^2. 
\ee 
Since $0<\mu <1$, one can choose  $M= M(c_0, |a|_\infty,  |b|_{\mu, \infty})>1$ large enough so that $\f{ 2^{2\mu }4 c_1}{M^{2(1-\mu)}}< \f{c_0} 4$. Consequently, the above estimate combined with  (\ref{lowerb1}) gives 
\be \label{lowerb2} 
|\w{u, \varphi_s^2 H_0 u} | \ge \f{c_0 }{4}(\|\nabla(\varphi_s u)\|^2 + \| \w{x}^{-\mu} \varphi_s u\|^2 ). 
\ee 
Remark that
\[
|\w{u, \varphi_s^2 H_0 u} | \le  \| \w{x}^{-\mu} \varphi_s u\| \| \w{x}^{\mu} \varphi_s H_0u\| \le \f{c_0}{8} \| \w{x}^{-\mu} \varphi_s u\| \|^2 +  \f{2}{c_0} \| \w{x}^{\mu} \varphi_s H_0u\|^2.
\]
It follows from  (\ref{lowerb2}) that \be
 \|\w{x}^{\mu} \varphi_s H_0 u \|^2 \ge \f{c_0^2}{16} (\| \w{x}^{-\mu} \varphi_s u\|^2
  +   \| \nabla (\varphi_s u)\|^2), \quad u\in C_0^\infty(\bR^n).
  \ee
By an argument of density, one obtains (\ref{e2.11}) with some constant $C>0$ independent of $s \in \bR$.
 \ef

\begin{cor} Under the conditions of Lemma \ref{lem2.2}, there exists some constant $C>0$ such that for
any $f \in L^{2, r}$ and   $u\in H^2_{\rm loc}$ such that $H_0u =f$, one has: $u \in L^{2, r-2\mu}$, $\nabla u \in L^{2, r-\mu}$ and \be \|\w{x}^{r-\mu}\nabla u\| + \|\w{x}^{r-2\mu} u\| \le C \|\w{x}^rf\|. \ee
\end{cor}
\pf It follows from Lemma \ref{lem2.2} with $s = \f{r-\mu}2$. \ef

Lemma \ref{lem2.2} shows that $H_0: D(H_0) \to R(H_0) :=  \mbox{\rm Range} (H_0) \subset L^2(\bR^n)$ is bijective. Let $G_0$ denote its algebraic inverse with $D(G_0) = R(H_0)$. Then one has 
\be 
H_0G_0 =1 \mbox{ on } R(H_0), \quad  G_0H_0 =1 \mbox{ on } D(H_0) 
\ee

\begin{lemma}  \label{lem3.1}
(a).  $G_0$ is a densely defined closed operator on $L^2(\bR^n)$. If $H_0$ is selfadjoint (resp., maximally dissipative),
 then $-G_0$ is also selfadjoint (resp., maximally dissipative). \\

 (b).  There exists some $C$ such that
 \be \label{e3.2}
\| \nabla(\varphi_s G_0 \varphi_{-s} \w{x}^{-\mu} w)\| + \|\w{x}^{-\mu}  \varphi_s G_0 \varphi_{-s} \w{x}^{-\mu} w)\| \le C\|w\| \ee for all 
$w \in \vD$  and $s\in \bR$. Here $\vD= \cap_{s\in \bR} L^{2,s}$.
\end{lemma}
\pf We firstly show that $D(G_0)$ is dense. Remark that $\Re H_0\ge 0$. Let $f \in \vD$ and $u_\ep = (H_0+ \ep)^{-1}f$, $\ep >0$.
 Since $\Re H_0 \ge 0$ and $H_0$ verifies the weighted coercive condition (\ref{ass2}), $H_0+\ep$ satsifies also (\ref{ass2}) with the same constant $c_0>9$ independent of $\ep>0$. Following the proof of Lemma \ref{lem2.2} with $H_0$ replaced  by $H_0+\ep$, one has that for any $s>0$
\[
\|\w{x}^{s-\mu}\nabla u_\ep\| + \|\w{x}^{s-2\mu} u_\ep\| \le C_s \|\w{x}^sf\|
\]
uniformly in $\ep>0$. For $s> 2 \mu$, this estimate implies that the sequence $\{u_\ep; \ep \in  ]0, 1]\}$ is relatively compact in $L^2$. Therefore there exists a subsequence $\{u_{\ep_k}; k \in \bN\}$ and $u \in L^2$ such that $\ep_k \to 0$ and $u_{\ep_k} \to u$ in $L^2$
 as $k \to +\infty$. It follows that $H_0u =f$ in the sense of distributions.
 The ellipticity of $H_0$ implies that $u \in H^2(\bR^n)$.
 Therefore $f \in R(H_0) = D(G_0)$. This shows that $\vD\subset D(G_0)$. In particular  $D(G_0)$ is dense in $L^{2,r}$ for any $r\in \bR$.
 The closeness of  $G_0$ follows from that of $H_0$. The other assertions can be easily checked.
\\

The argument used above shows that for any $w\in \vD$,  one can find $u\in D(H_0)$ such that $H_0 u = \varphi_{-s} \w{x}^{-\mu} w$. (\ref{e3.2}) follows from (\ref{e2.11}). \ef

Lemma \ref{lem3.1} shows that for any $s$, $\w{x}^{-\mu}  \varphi_s G_0 \varphi_{-s} \w{x}^{-\mu}$ defined on 
$\vD= \cap_{s\in \bR} L^{2,s}$ can be uniquely extended to a bounded operator on $L^2(\bR^n)$, or in other words, for any $s\in \bR$, $G_0 $ is bounded from $D(G_0) \cap L^{2,s}$ to $L^{2, s-2\mu}$:
\be
\|\w{x}^{-\mu}  \varphi_s G_0 u\| \le C \|\varphi_{s} \w{x}^{\mu} u\|
\ee
uniformly in $u \in D(G_0) \cap L^{2,s}$ and $s\in \bR$. This implies that $ G_0\vD \subset \vD$ 
 and $G_0$ extends to a continuous operator from $L^{2,s}$ to $L^{2, s-2\mu}$ for any $s\in\bR$. 
It follows that $ G_0^N(\vD) \subset \vD$ and by an induction, one can check that $G_0^N$ extends to a 
bounded operator from $L^{2,s}$ to $L^{2, s-2N\mu}$ for any $s \in \bR$. To simplify notation, we still denote $G_0$ (resp., $G_0^N$) its continuous extension by density as operator from $L^{2,s}$ to $L^{2,s-2\mu}$ (resp., from  $L^{2,s}$ to $L^{2,s-2N\mu}$). \\

\subsection{Gevrey estimates for the model resolvent}
\begin{theorem} \label{th3.2} Let $M>1$ be given in Lemma \ref{lem2.2}.   Denote
\be \label{scale}
x_{N,r} = \f x{R_{N,r}}  \mbox{ with }R_{N,r} = R_{(2N-1 + r)\mu} = M \w{ (2N-1 + r)\mu }^{\f 1 {1-\mu}}
\ee
where $N \in \bN$ and $r\in \bR_+$ and $M>0$ is a constant given by Lemma \ref{lem2.2}.  Set $\w{x_{N,r}} = (1 + |x_{N,r}|^2)^{\f 1 2}$. Then there exists some constant $C>0$ such that 
\be \label{e3.3}
 \| \w{x_{N,r}}^{-(2N+r)\mu} G_0^N  \w{x_{N,r}}^{ r\mu} \| \le C^N \w{ (2N-1 + r)\mu }^{\gamma N}, 
 \ee
 for any integer $N \ge 1$ and any $r\ge 0$. Here
 \be\label{gamma}
 \gamma = \f{2\mu}{1-\mu}.
 \ee
\end{theorem}
\pf  Making use of Lemma \ref{lem2.2}, one can check that operator
\be
 I_N =\w{x_{N,r}}^{-2N\mu - r\mu} G_0^N  \w{x}^{r \mu} 
 \ee
  is well defined on $\vD$ and  extends to a bounded operator on $L^2$. To show the estimate (\ref{e3.3}), we use an induction on $N$. Since $\w{x} \le \f 1 R\w{\f x R}$ for $R\ge 1$,   it follows from  (\ref{e3.2}) that
 \be \label{e3.4}
 \|\w{\f{x}{R_s}}^{-s-\mu}G_0 \w{\f{x}{R_s}}^{s-\mu}\|  \le C' R_s^{2\mu} \le C_1 \w{s}^{\gamma}
\ee uniformly in $s$, where $R = M\w{s}^{\f 1{1-\mu}}$. In particular, when $s = (1+r)\mu$, one has $ R_s = M\w{(1+r)\mu}^{\f 1{1-\mu}} =R_{1,r}$ and
 \be 
 \| I_1 \| \le C_1 \w{(1 + r)\mu}^{\gamma} 
 \ee 
 for all $r  \ge 0$, which proves (\ref{e3.3}) when  $N=1$. Assume now that $N\ge 2$ and  that 
 one has proved for some $C>0$  independent of $N$ and $r \ge 0$ that 
\be \label{e3.5} 
\| I_{N-1} \| \le C^{N-1} \w{(2N-3 + r)\mu}^{\gamma (N-1)}.
 \ee 
Write $I_N$ as
 \[
 I_N =  \w{x_{N,r}}^{-(2N+r)\mu} G_0  \w{x_{N-1,r}}^{(2N-2+r)\mu}\cdot I_{N-1} \cdot \w{x_{N-1,r}}^{-r \mu} \w{x_{N,r}}^{ r\mu}
 \]
 Notice that
\[
\w{x_{N,r}} \le \w{x_{N-1,r}} \le  \f{R_{N,r}}{R_{N-1,r}}\w{x_{N,r}}
\]
 for any $N \ge 2$. Applying (\ref{e3.4}) with $s = (2N-1+ r)\mu$, one obtains
 \[
 \| \w{x_{N,r}}^{-(2N+r)\mu} G_0  \w{x_{N,r}}^{(2N-2+r)\mu}\| \le C_1   \w{( 2N-1+r)\mu}^{\gamma}.
 \]
  Making use of the induction hypothesis, one can  estimate $I_N $ as follows:
\bea
\|I_N\| & \le &\| \w{x_{N,r}}^{-(2N+r)\mu} G_0  \w{x_{N-1}}^{(2N-2+r)\mu}\| \cdot \| I_{N-1}\|  \nonumber\\
& \le & \| \w{x_{N,r}}^{-(2N+r)\mu} G_0  \w{x_{N,r}}^{(2N-2+r)\mu}\|  \cdot \| (\f{\w{x_{N-1}}}{\w{x_{N,r}}})^{(2N-2+r)\mu}\| 
  \cdot \| I_{N-1}\|  \nonumber \\
&\le & C_1   \w{( 2N-1+r)\mu}^{\gamma} \cdot \left(\f{\w{(2N-1 +r)\mu}}{\w{(2N-3+r)\mu}}\right)^{\gamma (N-1 + \f r 2) } \cdot   C^{N-1} \w{( 2N-3 + r)\mu}^{\gamma(N-1)} \nonumber \\
&\le & C_1  \left(\f{2N-1 +r}{2N-3+r}\right)^{\gamma (N-1 + \f r 2) }   C^{N-1} 
 \w{( 2N-1+r)\mu}^{\gamma N}. \label{e3.6}
\eea 
The sequence $\{  \left(\f{2m-1 +r}{2m-3+r}\right)^{\gamma (m-1 + \f r 2) } ; m \ge 2 \}$ is uniformly bounded in $r \ge 0$. Hence there exists some $C_2>0$ such that
\[
 C_1  \left(\f{2m-1 +r}{2m-3+r}\right)^{\gamma (m-1 + \f r 2) } \le C_2
\]
for any $m \ge 2$ and $r\ge 0$. Increasing the constant $C$ if necessary, one can suppose without loss that 
$C_2 \le C$ and one obtains from (\ref{e3.6}) that 
\be 
\| I_N\| \le C^N \w{( 2N-1+r)\mu}^{N\gamma} 
\ee
 Theorem \ref{th3.2} is proven by induction.
\ef

Let $R_0(z)=(H_0-z)^{-1}$ denote the resolvent of $H_0$ and 
\[
 \Omega(\delta) =\{ z\in \bC^*; \f {\pi} 2 +\delta <  \arg z < \f {3\pi} 2 -\delta\} , \]
  $\delta >0$. Since $\Re H_0 \ge 0$, there exists some constant $C_1>0$ such that
\[
\|R_0(z)\| \le \f{C_1}{|z|}, \quad z\in \Omega(\delta).
\]
From the equation $R_0(z) =G_0 + zG_0 + z^2G_0^2R(z)$, it follows that as operators from $L^{2,s}$ to $L^{2,s-2\mu}$,  $s\in\bR$, one has
\be
\slim_{z\in \Omega(\delta), z\to 0} R_0(z) = G_0
\ee
for any $\delta >0$. Similarly one can check that for any $N \in \bN^*$, 
 one has
\be
\slim_{z\in \Omega(\delta), z\to 0} R_0(z)^N = G_0^N.
\ee
as operators from $L^{2,s}$ to $L^{2,s-2N\mu}$. By an abuse of notation, we denote $R(0) =G_0$. Thus  $R_0(z)$ is defined for $z$ in $\Omega(\delta)\cup\{0\}$. \\

\begin{cor} \label{cor3.3} The following Gevrey estimates of the resolvent hold.
\\

 (a). For any $a>0$, there exists some constant $C_a>0$ such that
 \be\label{e3.22}
 \| e^{-a\w{x}^{1-\mu}}  R_0(z)^N\| + \|   R_0(z)^N e^{-a\w{x}^{1-\mu}} \|\le  C_a^N N^{\gamma N}
 \ee
for any  integer $N\ge 1$ and $z \in \Omega(\delta)\cup\{0\}$.\\

(b).  Then there exists some constant $C>0$ such that
\be \label{e3.23} 
\| \chi(x)  R_0(z)^N\| + \|  R_0(z)^N\chi(x) \|\le C_\chi C^N N^{\gamma N}
 \ee 
 for any $\chi \in C_0^\infty(\bR^n)$,  $N\ge 1$ and $z \in \Omega(\delta)\cup\{0\}$. Here $\gamma = \f{2\mu}{1-\mu}$.
\end{cor}
\pf  Notice that $\|z R_0(z)\|$ is uniformly bounded in $\vL(L^2)$ for $z \in \Omega(\delta)$ ($\delta >0$ is  fixed) and that
\[
R_0(z)^{N} =  G_0^{N}(1+ z R_0(z))^N.
\]
According to Theorem \ref{th3.2} with $r=0$, one has for some constant $C>0$
\be \label{e3.25}
\| \w{x_{N,0}}^{-2N\mu} R_0(z)^N  \| \le C^N N^{\gamma N}, 
\ee
for any integer $N \ge 1$ and $z \in \Omega(\delta)\cup\{0\}$.
\\

Let $a>0$. Then 
\[
 \| e^{-a\w{x}^{1-\mu}}  R_0(z)^N\| \le \|e^{-a\w{x}^{1-\mu}}\w{x_{N,0}}^{2N\mu}\|_{L^\infty}  C^N N^{\gamma N}.
\]
To evaluate the norm $ \|e^{-a\w{x}^{1-\mu}}\w{x_{N,0}}^{2N\mu}\|_{L^\infty} $, consider the function
\[
f(r) = e^{-ar^{1-\mu}}\w{\f r{R_N}}^{2N\mu},
\]
where $r=|x|$ and $R_N = R_{N,0} = M\w{(2N-1)\mu}^{\f 1{1-\mu}}$.  One calculates:
\[
f'(r) = \f{f(r) }{r^\mu (R_N^2 + r^2)} ( -2a (1-\mu) (R_N^2 + r^2) + 2 N \mu r^{1+\mu}), r\ge 1.
\]
Let $A\ge 1$. Since $R_N \sim c'N^{\f 1{1-\mu}}$ for some constant $c'>0$, one can check that $ Nr^{1+\mu} \le \f{c}{A^{1-\mu}}r^2$ if $r\ge A R_N$ for some constant $c>0$ independent of $A$, $r$ and $N$. Therefore,  if $A= A(\mu, a)>1$ is chosen sufficiently large, one has
\[
f'(r) <0, \quad r > A R_N,
\]
thus $f(r)$ is decreasing in $[AR_N, +\infty[$. It is now clear  that
\[
\|e^{-a\w{x}^{1-\mu}}\w{x_{N,0}}^{2N\mu}\|_{L^\infty} \le \sup_{ 0 \le r AR_N} f(r) \le \w{A}^{2N\mu}
\]
This proves Part (a) of Corollary with $C_a = C \w{A}^{2\mu}$. \\

To prove  Part (b),  let $\chi \in C_0^\infty(\bR^n)$. Let $R>0$ such that supp $\chi \subset B(0,R)$. (\ref{e3.25}) shows that there exists some constant $C_1>0$ such that
\[
\|\chi(x) G_0^{N}(1+ z R_0(z))^N\| \le \|\w{x_{N,0}}^{2\mu N}  \chi  \|_{L^\infty}\times C_1^N N^{\gamma N}
\]
for any $\chi \in C_0^\infty(\bR^n)$,  $N\ge 1$ and $z \in \Omega(\ep)\cup\{0\}$.  Then One can check that
\[
\|  \w{x_{N,0}}^{2\mu N}  \chi\|_{L^\infty} \le \| \chi   \|_{L^\infty} (1+ \f{R^2}{M^2 ((2N-1)\mu)^{\f 2{1-\mu}}} )^{\mu N} \le C_2 2^{\mu N}
\]
for some constant $C_2$ depending only on $\chi$ and $R$, but independent of $N$. 
This proves (\ref{e3.23}) with   $C_\chi = C_2$ and $C=C_1 2^\mu$ which is independent of $\chi$.
\ef

Theorem \ref{th1.1} is a particular case of Corollary \ref{cor3.3}.  Corollary \ref{cor3.3} shows that the model resolvent $R_0(z)$ belongs to the Gevrey class of order $\sigma = 1+\gamma$ on $\Omega(\delta)\cup\{0\}$.

\sect{Quantum dynamics generated by the model operator}

\subsection{Subexponential time-decays of heat semigroup }

Consider now the model  operator of the form $ H_0= -\Delta + V_0(x)$  with 
$V_0(x) = V_1(x) - i V_2(x)$, $V_1(x), V_2(x)$ being real, satisfying Condition (\ref{ass2}). 
Denote $R_0(z)= (H_0-z)^{-1}$. Theorem \ref{th3.2} can be used to prove subexponential time-decay for local energies of solutions to the heat and   Schr\"odinger equations. To study the heat semigroup $e^{-tH_0}$, $t \ge 0$, we use Cauchy integral formula for semigroups and need some information of   the  resolvent on a contour in the right half complex plane passing through the origin. \\

\begin{prop} \label{prop4.1}
   Assume that $\Re H_0 \ge - a \Delta$ for some $a>0$ and
 that  the imaginary part of the potential $V_0(x)$ verifies the estimate
\be 
\label{ass1} 
|V_2(x)| \le C\w{x}^{-2\mu'}, \quad \forall x\in \bR^n, 
\ee
for some for some $0 <\mu'< \min\{\f n 2, 1\}$.
Then there exists some constant $C_0>0$ such that the numerical range $N(H_0)$ of $H_0$ is contained in
a region of the form $ \{z; \Re z \ge 0, |\Im z | \le C_0 (\Re z)^{\mu'}\}$. Consequently,  for any $A_0>C_0$
there exists some constant $M_0$ such that 
\be  \label{R4.1}
\|R_0(z)\| \le  \f{M_0}{|z|^{\f 1{\mu'}}} 
\ee 
for $z  \in \Omega :=\{z \in \bC^*; |z| \le 1, \Re z <0 \mbox{ or } \Re z \ge 0, |\Im z| >A_0 (\Re z)^{\mu'}\}$.
\end{prop}
\pf  
For $z = \w{u, H_0u} \in  N(H_0)$ where $u\in D(H_0)$ and $\|u\|=1$, one has 
\beas
\Re z & = & \Re  \w{u, H_0u} \ge a \|\nabla u\|^2 \\
|\Im z | & \le & \w{u, |V_2|u} \le C\|\w{x}^{-\mu'}u\|^2. 
\eeas 
According to the generalized Hardy inequality (\cite{hbst}), for $0<\mu' <\f n 2$ there exists some constant $C_{\mu'}$ such that
\be \label{Hardy} 
\|\w{x}^{-\mu'}u\|^2 \le \||x|^{-\mu'}u\|^2 \le C_{\mu'}
 \| |\nabla|^{\mu'}u\|^2.
\ee 
Let $\hat u$ denote the Fourier transform of $u$ normalized such that $\|\hat u\| = \| u\|$  and $\tau = \|\nabla u\|$. Then 
\beas
 \| |\nabla|^{\mu'}u\|^2& =& \||\xi|^{\mu'}\hat u\|^2
 = \||\xi|^{\mu'}\hat u\|_{L^2(|\xi|\ge \tau)} ^2 +  \||\xi|^{\mu'}\hat u\|_{L^2(|\xi|< \tau)} ^2 \\
&\le & \tau^{2(\mu'-1)} \||\xi| \hat u \|^2_{L^2(|\xi|\ge \tau)}  +  \tau^{2\mu'} \|\hat u\|_{L^2(|\xi|< \tau)} ^2 \\
& \le & 2 \tau^{2\mu'} = 2  \|\nabla u\|^{2\mu'}. 
\eeas  
This proves that $\Re z \ge 0$ and $|\Im z | \le C_0 (\Re z)^{\mu'}$ when $z \in N(H_0)$. The other assertions of Proposition are  immediate, since $\sigma(H_0) \subset \overline{N(H_0)}$ and 
\[ 
\|R_0(z)\| \le \f 1{\mbox{dist}(z, N(H_0))}.
\] 
\ef

Making use of the equation
\be
R_0(z)= \sum_{j=0}^{k-1} z^j G_0^{j+1} + z^k G_0^k R_0(z)
\ee
we deduce from Theorem\ref{th3.2} (with $r=0$ and $N=k$) and Proposition \ref{prop4.1} the following estimate for $k\in \bN^*$ and $k\ge \f 1{\mu'}$
\be
\|\w{x}^{-2k\mu }R_0(z)\| \le C
\ee
uniformly in $z\in \Omega$ and $z$ near $0$.
 Notice that under the conditions of Proposition \ref{prop4.1}, one can not exclude
 possible accumulation of complex eigenvalues towards zero. Making use of Proposition \ref{prop4.1}, one can prove  the following uniform Gevrey estimates in a domain located in the right half complex plane.\\

\begin{cor} \label{cor4.3}  Under the conditions of proposition \ref{prop4.1}, let $\kappa$ be an integer such that $\kappa + 1 \ge \f 1{\mu'}$. Then for any $a>0$ there exist $c_a, C_a>0$ such that
\be \label{3.31a} 
\| e^{-a\w{x}^{1-\mu}} \f{d^{N-1}}{dz^{N-1}} R_0(z)\| \le c_a C_a^N N^{(1+(1+\kappa)\gamma)N}, \quad \forall N \ge 1,
 \ee 
 and there exists some constant $C>0$ such that  for any  $\chi \in C_0^\infty(\bR^n)$, one has
\be \label{3.31} 
\| \chi(x) \f{d^{N-1}}{dz^{N-1}} R_0(z)\| \le C_\chi C^N N^{(1+(1+\kappa)\gamma)N}, \quad \forall N \ge 1,
 \ee 
uniformly in $z \in \Omega$. Here $\Omega$ is defined in Proposition \ref{prop4.1}.
\end{cor}
\pf  For $z \in \Omega$, decompose $R_0(z)$ into
\[
R_0(z) = A(z) + G_0^{\kappa +1} B(z)
\]
with $A(z) =\sum_{j=0}^\kappa z^j G_0^{j+1} $ and $B(z) = z^{\kappa +1}R_0(z)$.  By Proposition \ref{prop4.1}, 
$\|B(z)\|$ is uniformly bounded for $z\in \Omega$. Theorem \ref{th3.2} shows that for some constant $C_1$
\bea \label{3.32}
\| \w{x_{\kappa +1 ,r}}^{-(2\kappa +2 +r)\mu} G_0^{\kappa +1}  \w{x_{\kappa +1,r}}^{ r\mu} \| &\le& C_1 \w{ (2\kappa +1+ r)\mu }^{\gamma (\kappa +1)},  \\
\| \w{x_{\kappa +1 ,r}}^{-(2\kappa +2 +r)\mu} A(z) \w{x_{\kappa +1,r}}^{ r\mu} \| &\le & C_1 \w{ (2\kappa +1+ r)\mu }^{\gamma (\kappa +1)} \label{3.33}
 \eea
 for  any $r\ge 0$ and $|z| \le 1$.  Making use of the relation
 \[
 R_0(z)^N = A(z) R_0(z)^{N-1} +  G_0^{\kappa +1} R_0(z)^{N-1} B(z)
 \]
 one can show by an induction on $N$ that there exists some constant $C>0$ such that
 \be \label{3.34}
\| \w{x_{(\kappa +1)N, 0}}^{- 2(\kappa +1)N)\mu}R_0(z)^N \| \le C^{N} N^{N\gamma (1+\kappa)}
 \ee
 for any $N\ge 1$ and $z\in \Omega$. In fact, the case $N=1$ follows from (\ref{3.32}) and (\ref{3.33}).
 If (\ref{3.34}) is proven with $N$ replaced by $N-1$ for some $N \ge 2$, noticing that
 $x_{(\kappa +1)N, 0} = x_{\kappa +1, 2(\kappa+1) (N-1)}$, (\ref{3.32}) and (\ref{3.33}) with $r = 2(\kappa +1)(N-1)$ show that
\beas
\lefteqn{ \| \w{x_{(\kappa +1)N, 0}}^{- 2(\kappa +1)N)\mu}R_0(z)^N \|} \\
& \le & C_1 \w{ (2(\kappa +1)N-1)\mu }^{\gamma (\kappa +1)} ( \|\w{x_{(\kappa +1)(N-1),0}}^{ -2(\kappa +1)(N-1)\mu}R_0(z)^{N-1} \|   \\
& & +  \|\w{x_{(\kappa +1)(N-1),0}}^{ -2(\kappa +1)(N-1)\mu}R_0(z)^{N-1}\| \| B(z) \| ) \\
& \le & C_2 C^{N-1} N^{N\gamma (1+\kappa)}
\eeas
for some constant $C_2$ independent of $N$. Increasing the constant $C$ if necessary,  this proves (\ref{3.34}) for all $N \ge 1$ by an induction. (\ref{3.31a}) and (\ref{3.31}) are deduced from (\ref{3.34}) as in the proof of Corollary \ref{cor3.3}.
\ef

Note that in the applications,  we only use the Gevrey estimates at threshold.
As another consequence of Proposition  \ref{prop4.1}, we obtain the following estimate on the expansion of the resolvent at $0$:

\begin{cor} \label{cor4.2} Under the conditions of Proposition \ref{prop4.1}, assume in addition (\ref{ass2}) with  $\mu \in ]0,1[$. 
Then there exists some constant $c >0$ such that for any $z \in \Omega$ and $z$ near $0$, one has for some $N$ (depending on $z$) such that 
\be \label{3.45}
\| \w{x_{N,0}}^{-2N\mu}(R_0(z) -  \sum_{j=0}^N z^j G_0^{j+1}) \| \le  e^{-  c  |z|^{-\f 1 \gamma}}. 
\ee 
Here $\w{x_{N,0}}$ is defined in Theorem \ref{th3.2} with $r=0$.
\end{cor}
\pf Theorem \ref{th3.2} and Proposition \ref{prop4.1} show that for any $N$, one has 
\be \| \w{x_{N,0}}^{-2N\mu}(R_0(z) - \sum_{j=0}^N z^j G_0^{j+1}) \|
 \le C^N N^{\gamma N} |z|^{N+1  - \f 1 {\mu'}},
\ee 
for  all $z \in \Omega$ and $z$ near $0$. The remainder estimate can be minimized by choosing
an appropriate $N$ in terms of $|z|$.
For fixed $M'>1 $ and $z\neq 0$, take $N = [ \f 1{(CM'|z|)^{\f 1 \gamma}}]$. Then one has for $z$ in a small neighbourhood of zero and $z\neq 0$:
\[
C^N N^{\gamma N} |z|^{N+1  - \f 1 {\mu'}} \le e^{ - c_1N \log M }\le e^{ -c_2 |z|^{-\f 1 \gamma}}
\]
where $c_1, c_2$ are some positive constants. \ef

\begin{theorem} \label{th4.3}Let $H_0 =-\Delta + V_0(x)$ with $V_0\in \vV$.
Then for any $a>0$, there exist some constant $c_a, C_a>0$ suhc that
\be \label{eq3.11a}
\|e^{-a\w{x}^{1-\mu}}e^{-tH_0}\| + \|e^{-tH_0}e^{-a\w{x}^{1-\mu}}\| \le C_a e^{-c_a t^\beta}, \quad t>0,
\ee
with $\beta$ given by (\ref{beta}).
\end{theorem}
\pf Let $\Gamma$ be the contour  defined by $ \Gamma =\{z; \Re z \ge 0, |\Im z| = C (\Re z)^{\mu'}\}$ oriented in anti-clockwise sense, where $C>0$ is sufficiently large. Here $\mu'>0$ is appropriately small  such that both conditions (\ref{vV}) and (\ref{ass1}) are satisfied. By Proposition \ref{prop4.1}, the numerical range of $H_0$ is located on the right hand side of $\Gamma$ and one has 
\be 
e^{-tH_0} = \f i{2\pi} \int_{\Gamma} e^{-tz} R_0(z) dz. 
\ee 
Decompose $\Gamma $ as $\Gamma = \Gamma_0 + \Gamma_1$ where $\Gamma_0$ is the part of $\Gamma$ with $  0\le \Re z \le \delta$ while $\Gamma_1$ is the part of $\Gamma$ with $  \Re z > \delta$ where $\delta>0$ is sufficiently small. Clearly, the integral on $\Gamma_1$ is exponentially decreasing as $ t\to \infty$
\[
\| \int_{\Gamma_1} e^{-tz} R_0(z) dz \| \le Ce^{- \f t C}, t >0,
\]
for some constant $C>0$.
For $z \in \Gamma_0$, denote $f_N(z) = R_0(z) -  \sum_{j=0}^N z^j G_0^{j+1}$. Then 
\[
f_N(z) = z^{N+1} G_0^{N}R_0(z).
\] 
 Then Theorem \ref{th1.1} shows that for any $a>0$ there exist some constants $C, C_1>0$ such that
\be
\|e^{-a\w{x}^{1-\mu}}f_N(z)\| \le C_1 C^{N} |z|^{N+1-\f 1{\mu'}}N^{\gamma N}
\ee
for $z \in \Gamma_0$.  It follows that
 \beas
\lefteqn{\| \int_{\Gamma_0} e^{-tz} e^{-a\w{x}^{1-\mu}} R_0(z) dz \|} \\
& \le & \sum_{j=0} ^N \|e^{-a\w{x}^{1-\mu}} G_0^{j+ 1}\| | \int_{\Gamma_0} e^{-tz} z^j dz| +
\| \int_{\Gamma_0} e^{-tz} e^{-a\w{x}^{1-\mu}}f_N (z) dz \| \\
& \le &  C_2+  C_2\sum_{j=1} ^N  C^j j^{\gamma j} e^{-\delta t} + C_2 C^{N} N^{\gamma N} \int_{\Gamma_0}| e^{-tz}| |z|^{N+1 - \f 1{\mu'}} |dz | 
\eeas 
for some $C_2>0$ and for all $t>0$ and $N \ge 1$. Parameterizing $\Gamma_0$ by $z = \lambda \pm i c \lambda^{\f 1{\mu'}}$ with 
$\lambda \in ]0, \delta]$, one can evaluate the last integral as follows: 
\beas 
\int_{\Gamma_0}| e^{-tz}| |z|^{N+1 - \f 1{\mu'}} |dz | &\le&
C_3^N\int_0^\delta e^{-t \lambda} \lambda^{N+1 - \f 1{\mu'}} d\lambda \\
&\le & C_3^N  t^{-N -2 +\f 1{\mu'}}  \int_0^{\delta t} e^{-\tau} \tau^{N+1 - \f 1{\mu'}} d\tau \\
& \le &  C_4^N t^{-N -2 +\f 1{\mu'}} N^N 
\eeas 
for some $C_3, C_4 >0$. This proves that there exist some  constants $B_0$ and $B_1>0$ such that 
\be 
\| \int_{\Gamma_0} e^{-tz} \w{x_N}^{-2(N+1)\mu} R_0 (z) dz \| \le  B_0  B_1^N N^{\gamma N} (N e^{-\delta t} +   N^{ N} t^{-N -2 +\f 1{\mu'}}) 
\ee
 for any $t>0$ and $N \ge 1$. Choosing $N$ in terms of $t$ such that $N \simeq (\f t { M_1 B_1 })^{\f 1 {1+ \gamma}}$  as $t\to +\infty$ for  some fixed appropriate constant $M_1>1$, one obtains that 
\be 
\| \int_{\Gamma_0} e^{-tz}  \w{x_N}^{-2(N+1)\mu} R_0(z) dz \| \le Ce^{- \delta_0 t^{\f 1{1+\gamma}}}
 \ee 
 for some $C, \delta_0 >0$.  This proves that there exist some constants $C, c>0$   such that
 \be \label{e4.11}
 \|e^{-a\w{x}^{1-\mu}}e^{-t H_0}\|  \le C e^{- ct^{\beta }},  \quad  t>0,
 \ee
 with $\beta =  \f 1{1+\gamma} =\f{1+\mu}{1-\mu}$.  \ef
 
As a consequence of Theorem \ref{th4.3}, one obtains that  there exists some constant $c>0$ such that 
\be
 \|  e^{-t H_0}f\| \le C_R e^{- ct^{\beta }}\|f\|,  \quad t>0,
\ee
for all $f\in L^2(\bR^n)$  with support contained in $\{|x|\le R\}$, $R>0$.

\subsection{An estimate on spectral measure}

For the selafdjoint Schr\"odinger operator $H_0$ with a global positive and slowly decreasing potential $V_0$, it is known that under some additional conditions the spectral measure $E_0'(\lambda)$ of $H_0$ satisfies the estimate  that for any $N >0$ 
 \be
 E_0'(\lambda) = O_N(\lambda^N)
 \ee
 in appropriate spaces as $\lambda \to 0$ (see \cite{n}).  The Gevery estimates of the resolvent at threshold allow to give an improvement of this result. Let us begin with the following results on the boundary values of the resolvent up to real axis.
\\

\begin{lemma} \label{lem5.1}  Let $V_0(x) = V_1(x)- iV_2(x) $ with $V_1(x), V_2(x)$ real. Assume that  $V_1$ is of class $C^2$  on $\bR^n$ and that there exists
  $ \mu \in ]0, 1[$  and some constants $c_j >0$, $j=1,2,3$, such that
\bea
   c_1 \w{x}^{-2\mu}& \le &V_1(x) \le c_2\w{x}^{-2\mu},  \label{e3.16}  \\
   |(x\cdot\nabla)^j V_1 (x)| &\le & c_2\w{x}^{-2\mu}, \quad j = 1, 2  \label{e3.17} \\
 x\cdot\nabla V_1 (x) &\le & -c_3 \w{x}^{-2\mu}, \quad |x| > R  \mbox{  for some $R>0$},  \label{e3.18} \\
 |V_2(x)| &\le & c_2\w{x}^{-1-\mu-\ep_0}, \quad \ep_0 >0. \label{e3.19}
\eea
Then the eigenvalues of $H_0$ are absent in a neighbourhood of zero and
the boundary values of the resolvent $R_0(\lambda \pm i0) =\lim_{z\to \lambda, \pm \Im z >0} (H_0 -z)^{-1}$ exist for $\lambda \in [0, \delta]$ for some $\delta >0$ and are H\"older continuous as operators in $\vL(L^{2, \f{1+\mu} 2}; L^{2, -\f{1+\mu} 2})$.
\end{lemma}
\pf Let $H_{1} = -\Delta + V_1(x)$ be the selfadjoint part of $H_0$ and $R_{1}(z) =(H_{1} -z)^{-1}$. Then one knows  from \cite{n} that 
under the condition of this Lemma, $R_{1}(\lambda \pm i0)$ exists for $\lambda \in [0, \delta]$ for some $\delta >0$ and are Hölder continuous as operators in $\vL(L^{2, s}; L^{2, -s})$, $s>\f{1+\mu} 2$. 
Note that the smoothness assumption on the potential used in \cite{n} is only needed for higher order resolvent estimates.
\\

One knows that $ G_{0,1}=\lim_{z\to 0, z\not\in \bR_+}R_{1}(z)$ exists 
and  that $G_{0,1} V_2$ is a compact operator in $L^{2, -s}$ for $ \f{1 +\mu}{2} < s < \f{1 +\mu+\ep_0}{2}$.
 Therefore the kernel of $1+ i G_{0,1} V_2 $ is of finite dimension.
 From Lemma \ref{lem3.1} applied to $G_{01}$, one deduces that this kernel is contained in $L^{2,r}$ for any $r>0$.  Since $(1+ i G_{0,1} V_2)u=0$ if and only if $H_0u=0$, Lemma \ref{lem2.2} that $\ker (1+ i G_{0,1} v_2)$ in $L^{2, -s}$  is trivial.  Therefore $((1+ i G_{0,1} V_2)^{-1}$ is bounded in
$L^{2, -s}$. By the continuity of $R_{1}(z)$ for $z$ near $0$ and $z \not\in \bR_+$, one deduces that $1 +i R_{1}(z)V_2$ is invertible in $L^{2, -s}$ and its inverse is H\"older continuous  in $\vL(L^{2, -s})$ for $z$ near $0$ and $z \not\in \bR_+$. This implies in particular that the eigenvalues of $H_0$ are absent in a neighbourhood of zero and the limits $R_0(\lambda \pm i0) =\lim_{z\to \lambda, \pm \Im z >0} (H_0 -z)^{-1}$ exist for $\lambda \ge 0$ and small enough and are H\"older continuous in $\lambda$. \ef

\begin{cor} \label{th5.2}
Under the conditions of Lemma \ref{lem5.1}, assume in addition that  $V_2=0$ such that $H_0$ is selfadjoint. Denote by $E_0(\lambda)$ the spectral projection of $H_0$ associated with the interval $]-\infty, \lambda]$.  Let $s> \f{1+\mu}2$. Then  for any $a>0$, there exist some constants $c_a, C_a>0$ such that
\be \label{e3.20a} \|e^{-a\w{x}^{1-\mu}} E_0'(\lambda) \w{x}^{-s}\| \le C_ae^{-  c_a |\lambda|^{-\f 1 \gamma}}, \quad 0< \lambda \le \delta. 
\ee
\end{cor}
\pf Since $ E_0'(\lambda) = \f 1{2\pi i} ( R(\lambda+i0) - R(\lambda-i0))$,
$\| \w{x} ^{-s} E_0'(\la)\w{x}^{-s}\| $ is uniformly bounded   for $\lambda>0$ near $0$ (\cite{n}).
Iterating the resolvent equation, one obtains for any $N\in \bN$
\be
E_0'(\la) = \la^N G_0^N E_0'(\la), \quad 0 <\la \le \delta.
\ee
Applying (\ref{e3.3} with $r=s$, one deduces  as in the proof of Corollary \ref{cor3.3} that for any $a>0$, there exist some constants $c_a, C_a >0$ such that
\be
\|e^{-a\w{x}^{1-\mu}} E_0'(\lambda) \w{x}^{-s}\| \le C_a c_a^N N^{\gamma N} \la^N
\ee
for all $N\in \bN$ and $\la \in ]0, \delta]$. It remains to choose $N$ in terms of $\la>0$ (it suffices to take $N $ equal to integer part of 
$c \lambda^{-\f 1 \gamma}$ for some appropriate $c>0$) such that
\[
c_a^N N^{\gamma N} \la^N \le C'e^{-c'\la ^{-\f 1 \gamma}}, \quad  0 <\la \le \delta,
\] 
 for some constants $c', C'>0$. (\ref{e3.20a}) is proved.
 \ef

\subsection{Subexpontial time-decay of Sch\"odinger semigroup}

To obtain  subexponential time-decay of  solutions to the Schr\"odinger equation associated with $H_0$, we shall use both techniques of analytic dilation and analytic deformation to study quantum resonances of $H_0$. These different techniques give rise to the same set of quantum resonances.  See \cite{ac,hm,hun,si}.
\\

Firstly, we use the analytic dilation method to prove a global resolvent estimate in some sector below the positive real half-axis. For  $V_0\in \vA$, denote $\widetilde{H}_0(\theta) = -(1+\theta)^{-2} \Delta + V_0((1+\theta) x)$ for $\theta \in \bC$ and $\theta$ near $0$. Set $\widetilde{R}_0(z, \theta) = ( \widetilde{H}_0(\theta) -z)^{-1}$. For $\theta$ real, 
$\widetilde{R}_0(z, \theta)$ is holomorphic in $\bC_+$ and meromorphic in $\bC\setminus \bR_+$. 
Since $V_0 \in \vA$, there exists some constant $\delta>0$ such that $\{\widetilde{H}_0(\theta);  \theta \in \bC, |\theta| < \delta\}$ is a holomorphic family of type A. For $\Im \theta >0$
small enough, the resolvent $\widetilde{R}_0(z, \theta)$  defined for $z \in \bC_+$  with $\Im z >> 1$
can be meromorphically extended across the positive real half-axis $\bR_+$ into
the sector $\{ z; \arg z > -  \Im \theta\}$.
The poles of $\widetilde{R}_0(z, \theta)$ in this sector are by definition quantum  resonances of $H$ which are independent of $\theta$ (\cite{ac}).
\\

We begin with the following elementary Hardy type inequality. \\

\begin{lemma}\label{lem5.3} Let $n\ge 2$ and $0 <s<n-1$. One has
\be \label{ex3.23}
 \|\w{x}^{-1-\f s 2} u\|^2 \le \f 1{ 2 \sqrt{(n-1)(n-1-s)}} (\|\nabla u\|^2 + \|\w{x}^{-s}u\|^2) 
 \ee for all $u\in H^1(\bR^n)$.
\end{lemma}
\pf Let $x = r \omega$, $r\ge 0$ and $\omega\in \bS^{n-1}$. For $u \in \vS(\bR^n)$, denote
\[
 F(r, \omega) = \f{|u(r\omega)|^2 r^{n-1}}{\w{r}^s}.
\]
Then one has
\[
 F'_r(r, \omega) = \f{((n-1)(1+r^2) -s r^2)|u(r\omega)|^2 r^{n-2}}{\w{r}^{s +2} }+
 2 \f{ r^{n-1}}{\w{r}^s}\Re (u'_r(r\omega) \overline{u(r\omega)}).
\]
Here $F'_r(r, \omega)$ is the derivation of $F(r, \omega)$ with respect to $r$. Since for $n \ge 2$, one has
 $\int\int_{\bR_+\times \bS^{n-1}} F'(r,\omega) \; dr d\omega =0$, we deduce the identity
\be \int\int_{\bR^n} \f{ (n-1) +  (n-1-s)r^2}{r\w{r}^{s +2}}|u(x)|^2 \; dx
 = -2 \int\int_{\bR^n}  \w{x}^{-s} \Re (u'_r\overline u) \; dx
\ee for any $u \in \vS(\bR^n)$. Inequality (\ref{ex3.23}) follows from the trivial bounds
\[
  2\sqrt{(n-1)(n-1-s)} \le \f{ (n-1)  + (n-1-s)r^2}{r}
\]
and
\[
 -2 \int\int_{\bR^n}  \w{x}^{-s} \Re (u'_r\overline u) \; dx \le 2 \|u'_r\| \|\w{x}^{-s} u\| \le
 \|\nabla u\|^2+  \|\w{x}^{-s} u\|^2
\]
and an argument of density. \ef

\begin{lemma}\label{lem5.4} Let  $V_0 \in \vA$. 
Then there exists
some constant $c_0>0$ such that  for $\theta\in \bC$  with $|\theta|$ sufficiently small and $\Im \theta >0$, one has \be 
\sigma (\widetilde{H}_0(\theta))\cap \{z \in \bC; \Im z >0  \mbox{ or }   \arg z > - c_0 \Im \theta\} = \emptyset 
\ee 
and 
\be \label{e5.15}
\|\w{x}^{-2\mu}\widetilde{R}_0(z,\theta)\| \le \f{1}{c_0\Im\theta \w{z}} 
\ee 
for $z \in \bC$ with $\arg z > - c_0 \Im \theta$.
\end{lemma}
\pf We only consider the case $\theta = i \tau$ with $\tau>0$ small enough. Since $V_0= V_1-i V_2 \in \vA$, one has
\beas
V_0((1+i\tau)x) &= &V_1(x) + \tau x\cdot\nabla V_2(x) \\
 & & - i(V_2(x  )
 - \tau x\cdot\nabla V_1(x ) + O(\tau^2\w{x}^{-2\mu})
\eeas
for $\tau>0$ sufficiently small. Let $z = \w{u, \widetilde{H}_0(\theta) u}$, $u\in H^2$ with $\|u\|=1$. Then, 
\bea
 \Re z & = & \f{1-\tau^2}{(1+\tau^2)^2} \|\nabla u\|^2 + \w{u, (V_1(x) + O(\tau \w{x}^{-2\mu}))u },
   \\
\Im z &=&  -\f{2\tau}{(1+\tau^2)^2}   \|\nabla u\|^2 - \w{u, (V_2(x ) - \tau x\cdot\nabla V_1(x ) )u}
 \nonumber \\
  &  & +\w{u,  O(\tau^2\w{x}^{-2\mu})) u}.
\eea
 This implies that  there exists $c>0$ such that
\be
 \Re z \ge c (\|\nabla u\|^2 +\|\w{x}^{-\mu} u\|^2),
\ee 
for $\tau>0$ sufficiently small. If $R\in ]0, \infty]$, one has for some $c'>0$
\[
V_2(x ) - \tau x\cdot\nabla V_1(x ) \ge c' \tau \w{x}^{-2\mu},  \forall x\in \bR^n,
\]
which gives that 
\be
 \Im z \le - c'' \tau  (\|\nabla u\|^2 +\|\w{x}^{-\mu} u\|^2) 
\ee
 for some $c''>0$.  This shows that  $\Im z \le  - C \tau  \Re z$ ($ C = c''c^{-1}$) if $R\in ]0, +\infty]$.
 \\
 
 If $R=0$, one has $V_2 (x) \ge 0$ for all $x \in \bR^n$ and
\[
V_2(x) - \tau x\cdot\nabla V_1(x ) \ge c_3 \tau \f{x^2}{\w{x}^{2\mu +2}}, \quad \forall x\in \bR^n,
\]
for some $c_3>0$.  In this case, one has
 \be
 \Im z \le - C\tau (  \|\nabla u\|^2 + \w{u,  \f{x^2}{\w{x}^{2\mu +2}} u}) + C \tau^2 \|\w{x}^{-\mu}u\|^2.
\ee 
 Lemma \ref{lem5.3} with $s =\mu$ shows 
\[
 \f{1}{\w{x}^{2\mu +2}}  \le   \f 1{ 2 \sqrt{(n-1)(n-1-\mu)}}( -\Delta +  \f{1}{\w{x}^{2\mu}})
\]
in the sense of selfadjoint operators. For $0 <\mu  < \f 3 4$ when $n=2$ and $\mu \in ]0, 1[$ if $n \ge 3$, one has  
\[
\alpha \equiv  \f 1{ 2 \sqrt{(n-1)(n-1-\mu)}} <1.
\]
This proves that 
\beas 
\|\nabla u\|^2 + \w{u,  \f{x^2}{\w{x}^{2\mu +2}} u} 
& = & \|\nabla u\|^2 + \w{u,  (\f{1}{\w{x}^{2\mu }} - \f{1}{\w{x}^{2\mu +2}}) u}   \\
&\ge &
(1-\alpha) (\|\nabla u\|^2 + \|\w{x}^{-\mu } u\|^2). 
\eeas 
Consequently, one obtains 
\be
 \Im z \le  - C (1-\alpha) \tau  (  \|\nabla u\|^2 +  \|\w{x}^{-\mu } u\|^2) +
 C \tau^2 \|\w{x}^{-\mu}u\|^2 \le -C_1 \tau \Re z
\ee 
for some $C_1>0$ if $\tau >0 $ is small enough. 
This proves that the numerical range of $\widetilde{H}_0(\theta)$ is included in the region 
$ \Gamma= \{z;  \Re z \ge 0, \Im z \le  - C_1\tau \Re z \}$ for some constant $C_1>0$. 
Since $\sigma(\widetilde{H}_0(\theta)) \subset \Gamma$ and $\|\widetilde{R}_0(z, \theta)\| \le dist(z, \Gamma)^{-1}$, 
one has 
\be
\|\widetilde{R}_0(z, \theta)\| \le  \f{1}{c_0\Im\theta  \; |z|} 
\ee 
for $z \in \bC$ with $\arg z > - c_0 \Im \theta$ for some $c_0>0$.
The conclusion of Lemma \ref{lem5.4} follows now from Theorem \ref{th3.2}  applied  to $H_0= \widetilde{H}_0(\theta)$ which verifies Conditions (\ref{e1.1}) $-$ (\ref{ass2b}) uniformly in $\theta$ with $\Im \theta \ge 0$ and $|\theta|<\delta$. \ef

In order to obtain subexponential time-decay estimates for $\chi e^{-it H_0} \chi$,  $\chi \in C_0^\infty(\bR^n)$,  we use analytical distortion outside of $H_0$  the support of $\chi$. Let $R_0>1$  and $\rho \in C^\infty(\bR)$ with $ 0 \le \rho\le 1$ and $\rho(r) =0$ if $ r \le 1$ and 
 $\rho(r) =1$ if $ r \ge 2 $. Define for $R_0 >1$
 \be \label{5.16}
 F_\theta(x) = x(1 + \theta \rho(\f{|x|}{R_0})), \quad x \in \bR^n.
 \ee
When $\theta \in \bR$ with $|\theta|$ sufficiently small, $x\to F_\theta(x)$ is a global diffeomorphism on $\bR^n$.
Set
\be \label{5.17}
U_\theta f (x) = |DF_\theta(x)|^{\f 12} f(F_\theta(x)), \quad f \in L^2(\bR^n),
\ee
where $DF_\theta (x)$ is the Jacobi matrix and $ |DF_\theta(x)|$  the Jacobian of  the change of variables: $x \to F_\theta(x)$.
One has
\be
|DF_\theta(x)| = \left\{\begin{array}{ccc}
1, &  &  |x| < R_0; \\
  (1+\theta)^n, & \quad & |x|>2R_0 
\end{array} \right.
\ee
 $ U_\theta $ is unitary in $L^2(\bR^n)$ for $\theta $ real with $|\theta|$ sufficiently small.
Define the distorted operator $H_0(\theta)$ by
\be
H_0(\theta) = U_\theta H_0 U_\theta^{-1}.
\ee
One can calculate that
\be \label{e5.31}
H_0(\theta) = -\Delta_\theta + V(F_\theta(x))
\ee
where $-\Delta_\theta = ^t\nabla_\theta \cdot \nabla_\theta$ with
\be
\nabla_\theta = (^tDF_\theta)^{-1}\cdot \nabla  -\f 1 { |DF_\theta|^2}  (^tDF_\theta)^{-1}\cdot(\nabla |DF_\theta|)
\ee
In particular,  $ \nabla_\theta f = (1+ \theta)^{-1} \nabla  f$ if $f$ is supported outside the ball $B(0, 2 R_0)$. 
If $V_0 \in \vA$, $H_0(\theta)$ can be extended to a holomorphic family of type A  for $\theta$ in a small complex neighbourhood of zero. $H_0(\theta)$ and  $\widetilde{H}_0(\theta)$ coincide outside the ball $B(0, 2R_0)$ and they have the same essential spectra. In addition  their discrete  eigenvalues in the region  
$\{z\in \bC; \Re z \ge 0, \Im z >  - c\Im \theta \Re z\}$ for some $c>0$ small enough are the same (see \cite{ac, hun, hm}).
Since $\widetilde{R}_0(z, \theta)$ is holomorphic in $z$ there, so is $R_0(z, \theta) = (H_0(\theta) -z)^{-1}$.
\\

Remark that the distorted operator $H_0(\theta)$ satisfies the conditions (\ref{ass2}) and (\ref{ass2b}) with some constant $c_0>0$ independent of $R_0>1$. 
Consequently, Lemma \ref{lem2.2} applied to $H_0(\theta)$ implies that $\w{x}^{-2\mu} R_0(0, \theta)$ is defined on the range of $H_0(\theta)$  and extends to a bounded operator in $L^2(\bR^n)$ and Theorem \ref{th3.2} holds for $G_0(\theta) = R_0(0, \theta)$ for some constant $C$ independent of $ R_0$ used in the analytical distortion. \\

\begin{prop}\label{prop5.5} Assume the conditions of Lemma \ref{lem5.4}. Denote 
\[ 
\Omega_1(\theta)  = \{z\in \bC; \Re z > 0,   \Im z \ge -c   \Im \theta \Re z\}. 
\]
with $c>0$ appropriately small. Then $\Omega_1(\theta) $ is contained in resolvent set of $H_0(\theta)$ and there exists  some constant $C>0$
\be \label{e5.34}
\| \w{x}^{-2\mu} R_0(z, \theta) \w{x}^{-2\mu} \| \le \f C {\w{z}}, \quad z \in  \Omega_1(\theta).
\ee
\end{prop}
\pf  For $ z \in \Omega_1(\theta) $ and $|z| $ large, (\ref{e5.34}) follows from Lemma \ref{lem5.4} by  an argument of perturbation. 
 For $z \in \Omega_1(\theta) $ and $|z| $ small, one compares 
 $R_0(z, \theta)$ with $\widetilde{R}_0(z, \theta)$ and $ R_0(0, \theta)$.\\
 
 Let $\chi \in C_0^\infty(\bR^n)$ such that $ \chi(x) =1$ if $ |x| \le 2 R_0$. 
 On the support of $1 - \chi$, $ H_0(\theta) = \widetilde{H}_0(\theta)$. 
 For $z \in \Omega_1(\theta) $ and $|z|$ small, one has
 \beas
 R_0(z, \theta) & = &  R_0(0, \theta) + zR_0(0, \theta) R_0(z, \theta) \\
 & = &  R_0(0, \theta) + zR_0(0, \theta) ( \chi  (2-\chi) + (1-\chi)^2) R_0(z, \theta)\\
 & = &  R_0(0, \theta) + zR_0(0, \theta)  \chi  (2-\chi) R_0(z, \theta) \\
  &  &  +   zR_0(0, \theta) (1-\chi) \widetilde{R}_0(z, \theta) (1-\chi) \\
  & &  +   z R_0 (0, \theta) (1-\chi) \widetilde{R}_0(z, \theta) [ ( 1+\theta)^{-2} \Delta, \chi] R_0(z, \theta)
 \eeas
Recall that for $\Im \theta >0$, there exists some constant $C>0$ such that
\[
 \|\widetilde{R}_0(z, \theta) \w{x}^{-2\mu}\| \le   C, \mbox{ for } z \in \Omega_1(\theta).
\]
By the ellipticity of the operator, this implies that
\be
 \| \widetilde{R}_0(z, \theta) \nabla \w{x}^{-2\mu}\| \le   C, 
\ee
 for $ z \in \Omega_1(\theta)$ and $|z| \le 1$.
Therefore there exists possibly another constant $C$ such that
\be  \label{e5.35}
\|\w{x}^{-2\mu} R_0(z, \theta)\w{x}^{-2\mu}\|  \le  C + C |z| \|\w{x}^{-2\mu} R_0(z, \theta) \w{x}^{-2\mu} \| 
\ee
for $z\in \Omega_1(\theta)$ and $|z| \le 1$. This shows that 
$\|\w{x}^{-2\mu} R_0(z, \theta)\w{x}^{-2\mu}\|$ is uniformly bounded for  $z_0\in \Omega_1(\theta)$ and $|z| $ sufficiently small.
(\ref{e5.34}) is proved. \ef

\begin{theorem} \label{th5.6}  Let $V_0\in \vA$. There exists some constant $c>0$ such that
for any function $\chi\in C_0^\infty(\bR^n)$  there exists some constant $C_\chi>0$ such that
\be  \label{e5.60}
\|\chi(x)e^{-itH_0}\chi(x)\| \le C_\chi e^{-c |t|^{\beta}}, \quad t>0. 
\ee
\end{theorem}
\pf  Let $R_1>0$  such that $\supp \chi \subset B(0, R_1)$.
Let $U(\theta)$ be defined as before with $R_0 >R_1$. Then one has
\[
\chi(x)e^{-itH_0}\chi(x) = \chi(x)e^{-itH_0(\theta)}\chi(x)
\]
for $\theta \in \bR$ with $|\theta|$ small.  For $\theta \in \bC$  with $\theta$ near zero and $\Im \theta >0$, $H_0(\theta)$ is strictly sectorial and the resolvent $R(z, \theta)$ is holomorpic in $z\in \bC$ with
$ -c \Im \theta < \arg z  < \pi + c$ for some $c>0$. 
Making use of Proposition \ref{prop5.5}, one can check that
\be \label{e5.61}
\chi(x)e^{-itH_0}\chi(x) = \f i{2\pi} \int_{\Gamma'} e^{-it z}  \chi(x) R(z, \theta) \chi(x) dz
\ee
where
\[
\Gamma'= \{ z =  r e^{-i \delta  }; r \ge 0\}\cup \{z= - r e^{i\delta}, r\ge 0\}
\]
for $\delta = \delta(\Im \theta) >0$ small enough. $\Gamma'$ is oriented in anti-clockwise sense.
\\

The remaining part of the proof of (\ref{e5.60}) is the same as in Theorem \ref{th4.3} and will not be repeated here.
We just indicate that if one denotes $G_0(\theta) = R(0, \theta)$, then one has 
\[
R_0(z, \theta) = \sum_{j=0}^N z^j G_0(\theta)^j + z^{N+1} G_0(\theta)^N R_0(z, \theta)
\]
 for
$ z\in \Gamma'$ and $ z $ near $0$, and Theorem \ref{th3.2}  with $r=2$ and Proposition \ref{prop5.5} show that
\bea \label{5.27a}
\lefteqn{\|\chi(x) G_0(\theta)^N R(z, \theta)\chi(x)\|} \\
& \le  & C_\chi  \|\w{x_{N,0}}^{-2(N+1)\mu}G_0(\theta)^N R_0(z, \theta) \w{x_{N,0}}^{-2\mu}\| \le C_{\chi, \Im\theta} C^N N^\gamma \nonumber
\eea
with some constant $C$ independent of $\chi $ and $\theta$.  By choosing appropriately $N$ in terms of $t$ as in the proof of Theorem \ref{th4.3}, one obtain some $c>0$ independent of $\chi$ such that 
(\ref{e5.60}) holds.
\ef

Theorem \ref{th5.6} generalizes in particular   a result of D. Yafaev \cite{y1}  in one-dimensional selfadjoint  case to higher dimensions.\\

\sect{Compactly supported perturbations of the model operator}

 Consider non-selfadjoint Schr\"odinger operator $H$ of the form
\be\label{6.1}
H=H_0+W(x).
\ee
where $H_0 = -\Delta + V_0(x)$ with $V_0$ in $\vV$ or $\vA$ and $W\in L^\infty_{\rm comp}(\bR^n)$.
Then the essential spectrum of $H$ is equal to $[0, +\infty[$ and the possible accumulation points of  complex eigenvalues of $H$ are contained in  $\bR_+$.
We begin with the analysis of  positive resonances for a class of non-selfadjoint Schr\"odinger operators $H= -\Delta + V(x)$ with  $V(x)$ holomorphic outside some compact set. Since we are considering behavior of solutions as $t\to +\infty$, the main attention is paid to here outgoing positive resonances. Incoming positive resonances are invisible in the limit $t\to + \infty$.   
\\

\subsection{Positive resonances of non-selfadjoint Schr\"odinger operators}

Consider a class of non-selfadjoint Schr\"odinger operators  $H= -\Delta + V(x)$ which are
 compactly supported perturbations of $H_0 = -\Delta + V_0$ with
$\Im V_0(x) \le 0$ and $V_0(x)$ extends to  a holomorphic function a a region of the form
$\{x\in \bC^n; -x|>R, |\Im x| <\delta |Re x|\}$ and satisfies there
\be \label{6.1.0}
|V_0(x)| \le C\w{\Re x}^{-\rho}
\ee
for some constants $R, \delta, C, \rho >0$. Suppose in addition $(x\cdot\nabla_x)^jV_0 $, $j=0,1,2$, are $-\Delta$-compact.
Then the set $r_+(H)$ of outgoing positive resonances of $H$ is well defined by Definition \ref{def2.1}. 
For the model operator $H_0$, one has $r_+(H_0) =\emptyset$. Since $H_0$ is dissipative, the boundary value of thee resolvent
\be
R_0(\la +i0= \lim{z \to \la, \Im z>0} (H_0 -z)^{-1}
\ee
exists in $\vL(-1, s; 1,-s)$, $s > \f 1 2$, for $\la>0$  and is continuous in $\la$ (\cite{roy}). 
Let $U_\theta$ be the analytic distortion defined by (\ref{5.17}) with $R_0$ sufficiently large. Then $H_\theta$
\[
H_\theta = U_\theta H U_\theta^{-1}
\]
defined for $\theta$ real can be extended to a holomorphic family of type  A for $\theta$ in a complex neighbourhood of zero.
For $\Im \theta >0$, spectrum of $H_\theta$ is $\{z, \arg z > - c \Im \theta \}$ is discrete for some constant  $c>0$ and is independent of the function $\rho$ used in the distortion (\cite{hun}). 
\\

\begin{theorem}\label{th6.1.1} There exists some constant $\theta_0>0$ such that 
\be
\sigma_d(H(\theta))\cap \bR_+ = r_+(H)
\ee
for $ \Im \theta >0$ and $|\theta|< \theta_0$.
 In particular, outgoing positive resonances of $H$ is at most a countable set with zero as the only possible accumulation point.
\end{theorem}
\pf 
Let $\lambda_0>0$.  If $\lambda_0 \not\in r_+(H)$, then $-1$ is not an eigenvalue of the compact operator $R_0(\lambda_0 +i0)W$ on 
$L^{2, -s}$, $\f 1 2 < s < \rho -\f 1 2$. Therefore 
$(1+ R_0(\lambda_0 +i0)W)$ is invertible on $L^{2,-s}$.
Since $\lambda \to R_0(\la+i0) $ is continuous as operator from $L^{2, s}$ to $L^{2,-s}$, we deduce that $-1$ is not an eigenvalue of 
$R_0(\la +i0)W$ for $\lambda\in\bR$ and $\lambda$ sufficiently near $\lambda_0$. It follows that the boundary value of the resolvent
 $R(\lambda +i0) = \lim_{z \in \bC_+, z \to \lambda} R(z)$ exists   and
\be
\|\w{x}^{-s}R(\lambda+ i0 )\w{x}^{-s}\| \le C
\ee 
for $\lambda$ near $\lambda_0$ and $s> \f 1 2$. This proves that for any $\chi\in C_0^\infty(\bR^n)$, $\chi R(z)\chi$ is bounded for
 $z\in \bC_+$ and $z$ near $\lambda_0$. Therefore the meromorphic extension of $\chi R(z)\chi$ from the upper half complex plane is in fact holomorphic in neighbourhood of $\lambda_0$.  It follows that $ \sigma_d(H_\theta)\cap \bR_+\subset r_+(H)$. 
\\

 Conversely if $\lambda_0\in r_+(H)$, then there exists a non zero solution $u \in L^{2, -s}$ for any $s > \f 1 2$ such that
 \[
 u = - R_0(\la_0 +i0)Wu.
 \] 
 Let $R_0(z, \theta) = U_\theta R_0(z) U_\theta^{-1}$  and $u_\theta = U_\theta u$, for $\theta \in \bR$. Then one  has
 \be \label{6.1.2}
 u_\theta =   - R_0(\la_0 +i0, \theta)Wu
 \ee
  if the analytic distorsion is made outside the support of $W$. Since outgoing resonances
   of the dissipative operator $H_0$ are absent, $R_0(z, \theta)$ defined for $\Im z >0$ and $\theta$ real can be holomorphically extended for $\theta \in \bC$ with $\Im \theta >0$ and $|\theta| <\theta_0$ for some $\theta_0>0$ depending on domain of the analyticity of $V_0$. After this extension in $\theta$, $R_0(z, \theta)$ is holomorphic for $z$ near
 $\la_0$ and $\Im z > -c \Im \theta \Re z$ for some $c>0$. By (\ref{6.1.2}), $u_\theta$  can be extended in $\theta$  for $\Im \theta >0$.
 By (\ref{6.1.2}), $u_\theta \in L^2$  so long as $\Im \theta >0$. In addition $u_\theta \neq 0$ because
\be
\|\w{x}^{-s'}(u_\theta - u)\| \le |\w{x}^{-s}  (R_0(\la_0 +i0, \theta) -  R_0(\la_0 +i0, \theta)) Wu\| \le C |\theta|^\eta
\ee
for some $\eta>0$ if $s' > \f 1 2$. This proves $\la_0 $ is an eigenvalue of $H_\theta$ with $u_\theta$ as an eigenfunction when $\Im \theta >0$ and $|\theta|$ is small enough. Therefore $r_+(H) \subset \sigma_d(H_\theta) $. \ef

\begin{remark}\label{rmk6.0}
In Theorem \ref{th6.1.1}, the condition $\Im V_0 \le 0$ is used  because potential $V_0$ may  have   a long-range tail. Similarly if 
$\Im V_0 \ge 0$ one can  prove
\be
\sigma_d(H_\theta)\cap\bR_+ =r_-(H)
\ee
for $\Im \theta <0$ and $|\theta|<\theta_0$. In particular if $V_0(x)$ is  real when $x\in \bR^n$ and satisfies  (\ref{6.1.0}),  then real resonances of $H$ are at most a countable set with zero as the only possible accumulation point.
\end{remark}

The  sign restriction on $\Im V_0$ is not necessary if $V_0$ is of short-range. To be simple we only give a result in dilation analytic case. Let $H=-\Delta + V(x)$ with $V$  a dilation analytic short-range potential: $V_\theta(x) = V(e^\theta x)$ defined for $\theta$ real extends to a holomorphic function for $\theta$ in a complex neighbourhood of zero: 
\be \label{6.0.1}
|V_\theta(x)| \le \w{x}^{-\rho}, 
\ee
for $x\in \bR^n$ and $|\theta| <\theta_0$ for some $C, \theta_0>0$ and $\rho>1$. Then $H_\theta = -e^{-2\theta} \Delta + V_\theta$ is a holomorphic family of type A for $\theta \in \bC$, $|\theta|<\theta_0$

\begin{theorem} \label{th6.1.2} Let $V$ be  dilation analytic and of short-range. Then  for $|\theta|<\theta_0$ one has
\be \label{6.1.3}
\sigma_d(H_\theta)\cap \bR_+ = r_\pm(H), \pm \Im \theta >0.
\ee
In particular, positive resonances of $H$ are at most a countable set with zero as the only possible accumulation point.
\end{theorem}
\pf The inclusion $\sigma_d(H_\theta) \subset r_+(H),  \Im \theta >0$ can be proved in the same way as the first part of the proof of Theorem \ref{th6.1.1}.
It remains to show that $r_+(H) \subset\sigma_d(H_\theta),  \Im \theta >0$. Let $\la \in r_+(H)$. Then $-1$  is an eigenvalue of the compact operator $K = (-\Delta - \la-i0)^{-1}V$ on $L^{2,-s}$, $\f 1 2 < s <  \f \rho 2$. Let 
\[
K_\theta = (-e^{-2\theta}\Delta -\la -i0)^{-1}V_\theta =  e^{2\theta} (-\Delta -e^{2\theta}(\la +i0))^{-1}V_\theta,  \]
for $\Im \theta \ge 0$, $|\theta|<\theta_0$. $K_\theta$ is a family of compact operators on $L^{2, -s}$ continuous with respect to $\theta$ in the half disk
$D_+(0, \theta_0) =\{\theta \in \bC; \Im \theta \ge 0, |\theta |<\theta_0\}$ and holomorphic for $\theta$ in its interior. In addition
\be \label{6.0.4}
\|\w{x}^{-s} (K_\theta -K)\w{x}^{s}\| \le C |\theta|^\eta
\ee
for  $\f 1 2 < s <  \f \rho 2$ and for some $\eta>0$, because $(-\Delta - e^{2\theta}(\la + i0))^{-1}$ is H\"older continuous in 
$\theta \in  D_+(0, \theta_0)$ as operator-valued function from $L^{2,s}$ to $L^{2,-s}$. It follows that in any small neighbourhood of $-1$, $K_\theta$ has at least one eigenvalue $z_\theta$ for $\theta\in D_+(0, \delta)$ if $\delta>0$ is sufficiently small. Since $K_\theta$ and $K_{\theta'}$ are unitarily equivalent if $\Im \theta = \Im \theta'$, $z_\theta$ is
independent of $\Re \theta$. It follows that $z_\theta$ is independent of $\theta$ for $\Im \theta >0$ and $|\theta|<\delta$ (Theorem 1.9 in Ch. VII of  \cite{k}). Since $z_\theta \to -1$ as $\theta \to 0$, one concludes that $z_\theta =-1$ for $\theta \in  D_+(0 , \delta)$ if $\delta>0$ is small enough.
This proves that $\la$ is an eigenvalue of $H_\theta$ if $\Im \theta >0$ and $|\theta|<\delta$. Therefore $r_+(H) \subset\sigma_d(H_\theta)$, 
$\Im \theta >0$ which completes the proof of (\ref{6.1.3}) with sign $+$. The equality with sign $-$ can be proved in the same way. The last affirmation is immediate since the set of positive resonances of $H$ is equal to $r_+(H)\cup r_-(H)$. \ef 

In the above proof, we showed that if $\la\in r_+(H)$, then there exists $c>0$ such that $-1$ is the only eigenvalue of $K_\theta$ inside the disk $D(-1, c)$ for all $\theta \in D_+(\la, \delta)$. Then one can define the Riesz projection of eigenvalue $-1$ of $K_\theta$ by 
\be
\pi_\theta = \f 1{2\pi i} \int_{ |z+ 1| = \f c 2} (z-K_\theta)^{-1} dz, \; \theta \in D_+(\la, \delta).
\ee
 The following result is immediate. 

\begin{cor} \label{cor6.1.3} Assume the conditions of Theorem \ref{th6.1.2}. Let $\lambda>0$ be an outgoing resonance of $H=-\Delta + V$.
Denote $\pi_\theta$ the Riesz projection of eigenvalue $-1$ of $K_\theta$, $\theta \in D_+(\la, \delta)$, $\delta>0$.
Then as operators on $L^{2,-s}$, $\f 1 2< s < \f \rho 2$, $\pi_\theta $ is continuous for $\theta  \in D_+(\la, \delta)$ and holomorphic for $\theta$ in the interior of this half disk. 
\end{cor}

\subsection{Proof of Theorem \ref{th1.2}}

\begin{prop} \label{prop6.1} Let $H_0 = -\Delta + V_0(x)$ with $V_0\in \vV$. Let  $W\in L^\infty(\bR^n)$ with compact support and $H= H_0 + W(x)$. 
Assume that $0$ is not an eigenvalue of $H$. Then one has: \\

(a). There exist some constants $c_1, \mu'>0$ such that outside  the set
\[
\Omega_1 =\{z\in \bC; \Re z \ge 0 \mbox{ and } |\Im z | \le c_1 |\Re z|^{\mu'}\},
\]
 there are at most a finite number of discrete eigenvalues of $H$. There exists some $\delta>0$ such that
 \be \label{6.2a}
 \|R(z)\|\le \f C{|z|^{\f 1{\mu'}}}  \mfor z\not\in \Omega_1 \mand |z| <\delta.
 \ee 
 
 (b). The limit
\be
R(0) =\lim_{z\to 0, z\not\in \Omega_1} R(z)
\ee
exists in $\vL(L^{2, s}; L^{2, s-2\mu}$ for any $s\in \bR$ and  one has
\bea \label{e6.32}
\|e^{-a\w{x}^{1-\mu}} R(z)^N\| &\le &C_a^{N+1} N^{\gamma N},  \\
\|\chi R(z)^N \|&\le &  C_\chi  C^{N} N^{\gamma N} \label{e6.33}
\eea
for all $N \in \bN^*$ and $z\in \Omega_0 = \{z\in \bC; \Re z <0, |\Im z|< -M \Re z\} \cup \{0\}$, $M>0$. Here $a>0$ and 
$\chi\in C_0^\infty(\bR^n)$, $C_a, c_a, C_\chi$ are some positive constants and $C>0$ is  independent of $\chi$.
\end{prop}
\pf Note that $G_0 W$ is a compact operator and that $0$ is not an eigenvalue of $H$ if and only if $-1$ is not an eigenvalue of $G_0W$.
So if $0$ is not an eigenvalue of $H$,  operator $1+G_0W$ is invertible on $L^2$. From Proposition \ref{prop4.1}, one deduces that $1+R_0(z)W$ is invertible for $|z| $ small and $z \not\in \Omega_1$ This shows that $0$ is not an accumulation pint of 
$\sigma(H) \setminus \Omega_1$. In addition, $ z\to 1+ R_0(z)W$ is holomorphic in $\bC\setminus \Omega_1$. The analytic Fredholm Theorem
shows that $(1+ R_0(z)W)^{-1}$ is a meromorphic function with at most a discrete set of poles in $\bC\setminus \Omega_1$.
These poles are eigenvalues of $H$. Since  $0$ is not an accumulation point of eigenvalues of $H$ in $\bC\setminus \Omega_1$,  the number of eigenvalues of $H$ in $\bC\setminus \Omega_1$ is finite. 
(\ref{6.2a}) follows  from Proposition \ref{prop4.1} and the equation
\[
R(z) =  (1+R_0(z)W)^{-1}R_0(z).
\]

To prove that Gevrey estimates of the resolvent $R(z)$ at threshold, we  remark that if $F(z)$ and $G(z)$ are two bounded operator-valued functions on $\Omega_0$ satisfies the Gevrey estimates
\bea
\| F^{(N)}(z)\| &\le & A C_1^N (N!)^{\sigma} \label{6.34} \\
\| G^{(N)}(z)\| &\le & B C_2^N (N!)^{\sigma } \label{6.35}
\eea
for all $N \in \bN$ and $z \in \Omega_0$ and for  some $\sigma>1$ and $A, B, C_1, C_2 >0$, then $F(z)G(z)$  satisfies the Gevrey estimates
\be \label{6.36}
\| (FG)^{(N)}(z)\| \le  AB C_3^N (N!)^{\sigma } 
\ee
 for all $N \in \bN$ and $z \in \Omega_0$ where
 \be
 C_3 = D_\sigma \max\{C_1, C_2\}  \mbox{ with } D_\sigma =\sup_{N\in \bN} \sum_{j=0}^N \big (\f{j! (N-j)!}{N!}\big)^{\sigma -1} <\infty;
 \ee
  and if $F(z)$ is invertible for $z\in \Omega_0$ with uniformly bounded inverse: 
 \be
 \| F(z)^{-1}\| \le M
 \ee
uniformly in  $z\in\Omega_0$, then the inverse $H(z) = F(z)^{-1}$ satisfies the Gevrey estimates
\be
\| H^{(N)}(z)\| \le  M C_4^N (N!)^{\sigma } \label{6.37} 
\ee
 for all $N \in \bN$ and $z \in \Omega_0$, where $C_4 = M C_1 D_\gamma$.
  Denote  $G^\sigma(\Omega_0)$ the set of bounded operator-valued functions on $\Omega_0$ verifying Gevrey estimes of order $\sigma>1$. 
Since $e^{-a\w{x}^{1-\mu}}R_0(z)$ and $\chi R_0(z)$ belong to $G^\sigma$ with 
\be
\sigma = 1+\gamma =\f{1+\mu}{1-\mu},
\ee
 Estimates (\ref{e6.32}) and (\ref{e6.33}) follow respectively  from equations
\bea
e^{-a\w{x}^{1-\mu}} R(z) & = & (1 + e^{-a\w{x}^{1-\mu}}R_0(z) W e^{a\w{x}^{1-\mu}})^{-1} e^{-a\w{x}^{1-\mu}}R_0(z) \\
\chi R(z) & = & (1 + \chi R_0(z) W )^{-1}\chi R_0(z)
\eea
where $\chi\in C_0^\infty(\bR^n)$ is taken such that $\chi (x) =1$ on $\supp W$.
\ef

\begin{prop}\label{prop6.2} Let $V_0\in \vA$.
Assume that $0$ is not an eigenvalue of $H= H_0 + W(x)$. Then one has
\\

(a). There exists $\delta >0$ such that $H$ has at most  a finite number of eigenvalues in
\[
\Omega_\delta= \{z \in \bC\setminus \{0\};  -\delta \le \arg z \le \pi + \delta\}
\]
and for any $\chi\in  C_0^\infty(\bR^n)$, $ \chi R(z)\chi $ defined for $\Im z>0$ extends meromorphically into $\Omega_\d$ and there exists some constant $C_\chi, c >0$ such that
\be \label{6.2.1}
\| \chi R(z)\chi \| \le  C_\chi
\ee
for $z\in \Omega_\delta$ and $|z|<c$. \\

(b). The limit $R(0) = \lim_{z\in \Omega_\delta, z\to 0} R(z)$ in $\vL(L^{2, s+2\mu}; L^{2, s-2\mu})$ for any $s\in \bR$ and one has
\be  \label{e6.2.2}
\|\chi R(z)^N \chi\| \le   C_\chi  C^{N} N^{\gamma N} 
\ee
for any $N \in \bN^*$ and $z \in \Omega_\delta\cup \{0\}$ with $|z|<c$. 
\end{prop}
\pf Since $V_2 \ge 0$, one has for $\Im z>>1$
\be \label{6.2.3}
\chi R(z)\chi  =  (1 + \chi R_0(z) W )^{-1}\chi R_0(z)\chi
\ee
where $\chi\in C_0^\infty(\bR^n)$ is  such that $\chi (x) =1$ on $\supp W$. Let $U_\theta$ be defined by (\ref{5.17}) with $R_0>>1$ such that
supp $\chi \subset \{x; |x| <R_0\}$. Let $H_0(\theta) = U_\theta H_0 U_\theta^{-1}$ and  $R_0(z, \theta) = (H_0(\theta) -z)^{-1}$. 
Then one has for $\theta \in \bC$, $\Im \theta >0$ and $|\theta|$ small, 
\[
\chi R_0(z) W =  \chi R_0(z, \theta)W, \quad   \chi R_0(z)\chi =  \chi R_0(z, \theta)\chi.
\]
For a fixed $\theta\in \bC$ with $\Im \theta >0$,  Proposition \ref{prop5.5} shows that $\chi R_0(z) W $ and $\chi R_0(z) \chi$ are holomorphic in $\Omega_\delta$ for some $\delta>0$. The analytic Fredholm Theorem implies that $\chi R(z) \chi$ extends to  a meromorphic in $\Omega_\delta$ given by
\be \label{6.2.4}
\chi R(z)\chi  =  (1 + \chi R_0(z,\theta) W )^{-1}\chi R_0(z, \theta)\chi.
\ee
Zero is the only possible accumulation point of these poles. To show that zero is in fact not an accumulation point, we firstly prove that 
for each $\chi$, $-1$ is not an eigenvalue of the compact operator $\chi G_0 W$. In fact if $-1$ were an eigenvalue of $\chi G_0 W$, then $-1$ would also be an eigenvalue
of $WG_0\chi = J (\chi G_0 W)^*J$ where $J$ is the complex conjugaison: $Jf(x) = \overline{f(x)}$. Let $\psi \in L^2$ with $\psi \neq 0$ such that
\[
\psi = -WG_0\chi \psi.
\] 
Then $\chi \psi =\psi$, $G_0\psi \in L^2$  and $H G_0\psi = (1+ WG_0\chi)\psi =0$ which gives $\psi=\chi \psi =0$, since $0$ is not an eigenvalue of $H$ and $G_0$ is injective. This is in contradiction  with $\psi \neq 0$. therefore $-1$ is not an eigenvalue of $ \chi G_0 W$ and $1 + \chi G_0 W$ is invertible with bounded inverse. Secondly by an argument of compactness, one deduces that
if $\chi_R$, $R>1$,  is a family of  smooth cut-offs such that $\chi_R(x) =\chi(\f x R)$ for some function $\chi\in C_0^\infty(\bR^n)$ with $\chi(x) =1$ for $|x| \le 1$, then there exists some constant $C>0$ such that
\be
\|(1+\chi_RG_0W)^{-1} \|\le C
\ee
uniformly in $R>1$, because  $ \|(1+\chi_RG_0W)^{-1} -  (1+G_0W)^{-1} \| \to 0$ as $R \to +\infty$. 
According to Proposition \ref{prop5.5}, one has
\[
1+\chi_R R_0(z,\theta) W  = 1 +\chi_R G_0(\theta) W  +  O(|z|)
\]
in $\vL(L^{2,-2\mu}; L^{2, -2\mu}) $ for $z\in \Omega_\delta$, where $O(|z|)$ is uniform in $R>1$. Consequently there exists some constant $c>0$ independent of $R$ such that the inverse $(1+\chi_R R_0(z,\theta) W )^{-1}$ exists and is holomorphic for $z\in \Omega_\delta$ with $|z|<c$ and there exists some contant $C>0$ such that
\be \label{6.2.5}
\|(1+\chi_R R_0(z,\theta) W )^{-1}\|_{ \vL(L^{2,-2\mu}; L^{2, -2\mu})} \le C
\ee 
uniformly in $R>1$ and  $z\in \Omega_\delta$ with $|z|<c$.   This proves that there is no pole of the meromorphic extension of $\chi_R R(z)\chi$ in  $ \Omega_\delta \cap \{ |z| \le c\}$ for some $c>0$ independent of $R$.
Since discret eigenvalues  and positive outgoing resonances of $H$  are poles of some cut-off resolvent $\chi_R R(z)\chi_R $ if $R>1$ is large enough, (\ref{6.2.5}) implies
\be
(\sigma_d(H) \cup r_+(H)) \cap \{z; z \in \Omega_\delta, |z| \le c \} =\emptyset
\ee
This proves the finiteness of eigenvalues of $H$ in $\Omega_\delta$, because zero is  the only possible accumulation point of eigenvalues of $H$  in $\Omega_\delta$. Estimate (\ref{6.2.1}) follows from (\ref{e5.34}) and (\ref{6.2.4}).
\\

Part (b) can be derived from (\ref{6.2.4}), Proposition \ref{prop5.5} and Theorem \ref{th1.1} applied to $G_0(\theta)$.
\ef

Since for $V_0 \in \vA$, $\Im V_0 \le  0$, Theorem \ref{th6.1.1} can be applied which implies that  zero is the only possible accumulation  point of $r_+(H)$. 
As consequence of Theorem \ref{th6.1.1} and Proposition \ref{prop6.2}, one obtains the following

\begin{cor} \label{cor6.1.2}  Assume the conditions of Proposition \ref{prop6.2}.  Then the set $r_+(H)$ is finite.
\end{cor}

{\noindent \bf Proof of Theorem  \ref{th1.2}} (a).  Theorem \ref{th1.2} (a) can be proved in the same way as  Theorem \ref{th4.3} for the model operator $H_0$. By Proposition \ref{prop6.1}, one can find a contour $\Gamma $ in the right half complex plane of the form
\[
\Gamma =\{z; \Re z \ge 0, |\Im z| = C (\Re z)^{\mu'}\}
\]
for some $C, \mu'>0$ such that $\sigma(H) \cap \Gamma =\{0\}$ and there are only a finite number of complex eigenvalues of $H$ located on the left of $\Gamma$. Let 
\[
\Lambda = \sigma(H) \cap  \{z; \Re z < 0 \mor \Re z \ge 0 \mand |\Im z| >C (\Re z)^{\mu'}\}. 
\]
Then one has
\be 
e^{-tH} -\sum_{\lambda \in \Lambda} e^{-tH}\Pi_\lambda = \f i{2\pi} \int_{\Gamma} e^{-tz} R(z) dz. 
\ee 
where $\Pi_\lambda$ is the Riesz projection of $H$ associated the eigenvalue $\lambda$. 
Making  use of Proposition \ref{prop6.1}, one can prove as in Theorem \ref{th4.3} that 
\be
\| e^{-a\w{x}^{1-\mu}} (e^{-tH} -\sum_{\lambda \in \Lambda} e^{-tH}\Pi_\lambda)\| \le C_a e^{-c_a t^\beta}, \quad t>0.
\ee
(\ref{eq1.8}) follows since if $\lambda  \in \sigma_d(H)$ with $\Re \lambda >0$, $\|e^{-tH} \Pi_\lambda\|$ decreases exponentially.
(\ref{eq1.9}) is deduced in a similar way. 
\\

 (b). According to Proposition \ref{prop6.2} and Corollary \ref{cor6.1.2}, there exists some $ \eta >0$ such that
$\Omega_\eta$ contains no poles with negative imaginary part of meromorphic extension  of $\chi R(z)\chi$ from the upper half complex plane.  Under the assumptions Theorem \ref{th1.2} (b), $\chi R(z)\chi$ has only a finite number of poles in $\overline{\bC}_+$ which are either discrete eigenvalues or positive outgoing resonances of $H$.
 Making use of (\ref{6.2.4}) for some fixed $\theta \in \bC_+$, one obtains the the representation formula:
\be
\chi e^{-itH}\chi  -\sum_{\lambda \mbox{ poles in } \overline{\bC}_+} {\rm Res } (e^{-itz}\chi R(z)\chi; \lambda) = \f i {2\pi } \int_{\Gamma_\eta}e^{-itz} \chi R(z)\chi dz, \quad t >0,
\ee
where  ${\rm Res } (e^{-itz}\chi R(z)\chi; \lambda)$ is the residue of $e^{-itz}\chi R(z)\chi$ at pole $\lambda$ and $\Gamma_\eta =\{ z = re^{-i\eta}; r \ge 0\} \cup \{ z = -re^{i\eta}; r \ge 0\}$ where $\eta>0$ is chosen such that there are no eigenvalues of $H$ on $\Gamma_\eta$, nor  between  $\Gamma_\eta$ and the real axis. 
It is easy to see that
\[
{\rm Res } (e^{-itz}\chi R(z)\chi; \lambda) = 
\chi e^{-itH}\Pi_\la \chi,   \mbox{ if } \lambda \in \sigma_d(H),
\]
while for $\lambda \in r_+(H)$ we can only affirm that
\[
{\rm Res } (e^{-itz}\chi R(z)\chi; \lambda) = 
\chi e^{-it\lambda}P_\la(t) \chi
\]
for some operator of finite rank $P_\la(t)$ which is polynomial in $t$. See Remark \ref{rmk7.1} on the rank of
coefficients of $P_\la(t)$.
The subexpoential time-decay of the contour integral can be deduced from the Gevrey estimates for the cut-off resolvent  (\ref{e6.2.2}).  
The details are the same as in the proof of Theorem \ref{th5.6} and are omitted here.
\ef

\subsection{Proof of Theorem \ref{th1.3}}

Assume now that zero is an eigenvalue of $H= H_0 + W$. Then $-1$ is an eigenvalue of $G_0W$ and $\ker (1 + G_0 W)$ in $L^{2}$ coincides with the eigenspace of  $H$ with eigenvalue zero.  We begin with a decay estimate of the associated eigenfunctions.\\

\begin{lemma}\label{lem6.3} Assume that $H_0 = -\Delta + V_0$ satisfies the condition (\ref{ass2}). Then there exits some constant $\alpha_0>0$ such that if 
$u\in H^2(\bR^n)$ such that $Hu =0$, then $e^{\alpha_0 \w{x}^{1-\mu}}u \in L^2(\bR^n)$.
\end{lemma}
\pf Let $\varphi(x) =\alpha \w{x}^{1-\mu}$, $\alpha >0$. Let $\chi$ be a smooth cut-off on $\bR^n$  such that $0 \le \chi \le 1$, $\chi(x)=1$ for  $|x| \le 1$ and $\chi(x)=0$ for  $|x| \ge 2$. Set 
\[
\varphi_R(x)= \chi(\f x R) \varphi(x), \quad R \ge 1.
\]
Then
\[
|\nabla \varphi_R(x)|  \le \alpha (1 +  \f{C}{R^{1-\mu}}) \w{x}^{-\mu} \le 2 \alpha \w{x}^{-\mu}
\]
 uniformly in $R \ge R_1$ where $R_1$ is sufficiently large. Let  $u \in H^2(\bR^n)$. Then  one  has
\bea
|\w{e^{2\varphi_R} H u, u}| &= &|\w{H(e^{\varphi_R} u), e^{\varphi_R} u} -\w{[\Delta, e^{\varphi_R}]u,  e^{\varphi_R} u} |\\
  &= & | \w{H(e^{\varphi_R} u), e^{\varphi_R} u} +2\w{ ( |\nabla\varphi_R|^2 e^{\varphi_R}u - \nabla\varphi_R\cdot\nabla (e^{\varphi_R}u),  e^{\varphi_R} u} | \nonumber \\
& \ge & |\w{H(e^{\varphi_R} u), e^{\varphi_R} u}| -(8\alpha^2 +2 \alpha) \|\w{x}^{-\mu} e^{\varphi_R}u \|^2 -  2\alpha  \|\nabla (e^{\varphi_R}u)\|^2 \nonumber
 \eea
uniformly in $R \ge R_1$.
Since $W$ is compactly supported,$\varphi_R$ is bounded on supp $W$ uniformly with respect to $R$.  Making use of the condition (\ref{ass2}), one obtains for some constants $c_0, C>0$, 
\be
|\w{H(e^{\varphi_R} u), e^{\varphi_R} u}| \ge c_0( \|\nabla (e^{\varphi_R}u)\|^2 + \|\w{x}^{-\mu} e^{\varphi_R}u \|^2) -C \|u\|^2
\ee
 for all $R\ge R_1$. For  $\alpha>0$ appropriately small, one deduces that
 there exists some constant $C_1>0$ such that
\be
\|\w{x}^{-\mu} e^{\varphi_R}u \|^2 +\|\nabla (e^{\varphi_R}u)\|^2 \le C_1 ( |\w{e^{2\varphi_R} H u, u}| + \|u\|^2)
\ee
for any $u \in H^2$ and $R \ge R_1$. If $u \in H^2$ such that $Hu=0$, it follows that
\be
\|\w{x}^{-\mu} e^{\varphi_R}u \|^2 +\|\nabla (e^{\varphi_R}u)\|^2 \le C_1  \|u\|^2
\ee
for all  $R \ge R_1$. This proves that
 $\w{x}^{-\mu}e^\varphi u  \in L^2(\bR^n)$ and $ \nabla (e^{\varphi}u)  \in L^2$. Lemma \ref{lem6.3} is proved, provided that $0<\alpha_0<\alpha$.
\ef

\begin{theorem}\label{th6.4} Let $H_0 = -\Delta + V_0(x)$ and $H = H_0 + W(x)$
with $V_0\in \vV$ and $W\in L^\infty_{\rm comp}$. Assume that $0$ is an eigenvalue of $H$ and that  both $H_0$ and $H$ are selfadjoint.
Let $\Pi_0$ denote the eigenprojection of $H$ associated with eigenvalue zero. Then there exist some constants $C,\mu', \delta >0$ such that 
\be \label{6.5}
R(z) = -\f{\Pi_0}{z} + R_1(z)
\ee  
for $z\in \Omega_1(\delta)$ where
\be 
\Omega_1(\delta) =  \{z\in \bC; |z| <\delta, \mbox{ either } \Re z <0 \mor \Re z \ge 0 \mand |\Im z | > C|\Re z|^{\mu'}\}.
\ee
The remainder $R_1(z)$ satisfies the estimates
\be \label{6.3.2}
\|\w{x}^{-s} R_1(z)\|  + \|R_1(z) \w{x}^{-s} \| \le C_s
\ee
for $s> 2\mu + \f 1 {\mu'}$ and $z\in \Omega_1(\delta)$; and  for any  $ a, M >0$ there exist some constants $C_a, c_a >0$ such that
\be \label{6.6}
\| e^{-a\w{x}^{1-\mu}} R_1^{(N)}(z)\| + \|   R_1^{(N)}(z)e^{-a\w{x}^{1-\mu} }\| \le C_a c_a^N N^{\sigma N},
\ee
for any $N \in \bN$ and $z \in \Omega_-$ where  $\Omega_- =\{ z; \Re z <0 \mand |\Im z |\le  -M \Re z\} \cup\{0\}$, $M>0$.  Here 
and in the following, 
$R_1^{(N)}(z)$ denotes the $N$-th derivative of $R_1(z)$ and $\sigma = 1 +\gamma =\f{1+ \mu}{1-\mu}$.
\end{theorem}
\pf We use the Grushin method to study the low-energy asymptotics for the resolvent of $H$
 by using the  equation
\be \label{R}
R(z)  = (1+ R_0(z) W)^{-1}R_0(z).
\ee
Since the method is well-known in selfadjoint case, we shall skip over some details and emphasize on the Gevrey estimates of the remainder.
Note that 
$  \ker_{L^{2,s}} (1+  G_0 W)  $ is independent  of  $s\in \bR$ and coincides with the eigenspace of $H$ associated withe the eigenvalue $0$. We need only to work in $L^2(\bR^n)$. \\

Let $\psi_1, \cdots, \psi_m$ be a basis of  $\ker(1+  G_0 W)$  such that 
\be \label{psij}
\w{\psi_j, - W \psi_k} =\delta_{jk}, \quad j, k =1, \cdots, m.
\ee
(\ref{psij}) can be realized because the quadratic form $\phi \to \w{\phi, -W \phi} =\w{\phi, H_0\phi}$ is positive definite on 
$\ker ( 1+ G_0 W )$. Define $Q : L^{2} \to L^{2}$ by
\be
Q f = \sum_{j=1}^m \w{-W \psi_j, f } \psi_j, \quad f \in   L^{2}.
\ee
Set $Q' =1-Q$. Then $Q$ commutes with $1+G_0W$. $-1$ is not eigenvalue of compact operator $Q'(G_0W)Q'$, hence $ Q'(1+G_0W)Q'$ is  invertible on the range of $Q'$ with bounded inverse. From Theorem \ref{th3.2} with $N=1$ and Proposition \ref{prop4.1}, one deduces that
\be
(R_0(z) -G_0)W = O(|z|)
\ee
for $z\in \Omega_1(\delta)$. It follows that if $\delta>0$ is small enough, 
\be
E(z) = (Q'(1+  R_0(z) W)Q')^{-1}Q'
\ee
is well-defined and continuous in $z \in \Omega_1(\delta)$ and is uniformly bounded:
\be
\|E(z)\| \le C
\ee
uniformly in  $z \in \Omega_1(\delta)$. By Corollary \ref{cor3.3} and (\ref{6.36}), $E(z)$ satisfies  Gevrey estimates
\be\label{6.3.1}
\| E^{(N)}(z)  \| \le C C'^N N^{\sigma N}
\ee
for some $C'>0$ and for all $z \in \Omega_-$. 

Define   $S : \bC^m \to D(H)$ and $T : L^2 \to \bC^m$ by 
\begin{eqnarray*}
S   c &=& \sum_{j=1}^m c_j\psi_j, \quad  c = (c_1,
                  \cdots, c_m)\in \bC^m, \\
Tf &= &(\w{-W\psi_1, f},  \cdots, \w{-W\psi_m, f}) \in \bC^m,\quad
f\in L^2.
\end{eqnarray*}
Set $W(z) = (1+  R_0(z) W)$ and 
\begin{eqnarray}
 E_+(z) &=& S - E(z) W(z)S, \\
E_-(z )  &=&  T - TW(z) E(z),\\
  E_{-+}(z)& =& - TW(z)S + T W(z) E(z) W(z)S.
\end{eqnarray}
Then one has the formula
\begin{equation} \label{W}
(1+ R_0(z) W)^{-1} = E(z) - E_+(z) E_{-+}(z)^{-1} E_-(z) \mbox{ on } H^{1, -s}.
\end{equation}
Since $E(z)$, $W(z)$ satisfy Gevrey estimates  of the form (\ref{6.3.1}) on $\Omega_-$, 
$E_\pm(z)$ and $E_{-+}(z)$ satisfy similar Gevrey estimates on $\Omega_-$. The leading term of $E_{-+}(z)$ can be explicitly calculated:
\be
E_{-+}(z) = -z \Psi + z^2 r_1(z)
\ee
where   the matrix $\Psi= \left(\w{\psi_j, \psi_k} \right )_{1\le j, k \le m}$  is positive definite and $r_1(z)$ satisfies the Gevrey estimates in $\Omega_-$. Consequently,
\be
E_{-+}(z)^{-1} = -\f{\Psi^{-1}}{z} + \widetilde{r}_1(z)
\ee
with  $\widetilde{r}_1(z)$ uniformly bounded on $\Omega_1(\delta)$  and $\widetilde{r}_1(z)$ satisfying the Gevrey estimates of the form (\ref{6.3.1}) in $\Omega_-$. Consequently $(1+ R_0(z) W)^{-1}$ is of the form
\be
(1+ R_0(z) W)^{-1} = \f{A_0}{z} + B(z)
\ee
where 
\be
A_0 = S\Psi^{-1} T
\ee
 is an operator of rank $m$ and $B(z)$ is uniformly bounded $\Omega_1(\delta)$ and satisfies 
the Gevrey estimates 
\be \label{6.3.3}
\| B^{(N)}(z) \| \le C C'^N N^{\sigma N}, \quad \forall N \in \bN, 
\ee
for $z$ in $\Omega_-$. From the equation $R(z) = (1 + R_0(z) W)^{-1}R_0(z)$ and Corollary \ref{cor3.3},  we deduce that
\be
R(z)  = -\f{\Pi_0}{z} + R_1(z)
\ee  
where $R_1(z)$ satisfies 
\[
\| R_1(z) \w{x}^{-2k\mu} \|\le C
\]
for $z\in \Omega_- $ if $k \in \bN$ and $k \ge \f 1 {\mu'}$ and 
\[
\|   R_1^{(N)}(z) e^{- a\w{x}^{1-\mu}}\| \le C_a c_a^N N^{\sigma N}
\]
for $z$ in $\Omega_-$ and $N \in \bN$. This proves  (\ref{6.3.2}) and (\ref{6.6}).
\ef

\begin{theorem}\label{th6.4b} Let $H_0 = -\Delta + V_0(x)$ and $H = H_0 + W(x)$
with $V_0\in \vA$ and $W\in L^\infty_{\rm comp}$. Assume that $0$ is an eigenvalue of $H$ and that  both $H_0$ and $H$ are selfadjoint.
Let $\Pi_0$ denote the eigenprojection of $H$ associated with eigenvalue zero. Let $\Omega_\delta$ be defined as in Proposition \ref{prop6.2}
and
 $\Omega_\delta(c) =\Omega_\delta\cap \{|z| <c\}$. Then there exist some constants $C, c,\mu', \delta >0$ such that for any 
 $\chi\in C_0^\infty(\bR^n)$ the cut-off resolvent $\chi R(z)\chi$ defined for $\Im z>0$ extends to a holomorphic function in 
 $\Omega_\delta(c)$ and one has
\be \label{6.5}
\chi R(z)\chi = -\f{\chi \Pi_0 \chi}{z} + R_2(z)
\ee  
for $z\in \Omega_\delta(c)$ where the remainder $R_2(z)$ is continuous up to $z=0$ and satisfies the Gevrey estimates
\be \label{6.6b}
\|\chi R_2^{(N)}(z) \chi\| \le C_\chi C^N N^{\sigma N}
\ee
for $z \in \Omega_\delta(c) \cup\{0\}$ and for all $N \in \bN^*$.  
\end{theorem}
\pf It suffices to prove (\ref{6.5}) for $\chi\in C_@\infty(\bR^n)$ with sufficiently large support.
 Let  $\chi_0 \in C_0^\infty(\bR^n)$ such that $0 \le \chi_1(x) \le 1$, $\chi_1(x) =1$ for $|x|\le 1$ and $0$ for $|x| \ge 2$. 
Set 
\be
\chi_j (x) = \chi_0(\f x{jR}), \quad j =1,2,
\ee
where $R >R_0$ is to be adjusted  and  $R_0$ is such that supp $W\subset \{x; |x| \le R_0\}$. Then $\chi_j W = W$ and $\chi_1 \chi_2 =\chi_1$. 
Then one has
\be  \label{6.61}
\chi_1 R(z) \chi_1 =(1+\chi_2 R_0(z, \theta)W)^{-1} \chi _2R_0(z, \theta) \chi_1,
\ee
where the analytical distortion is carried out outside the support of $\chi_2$. (\ref{6.61}) initially valid for $\theta$  real and $\Im z >0$ allows to extend $z \to \chi R(z) \chi$ into a sector below the positive real axis when $\Im \theta >0$. In the following  $\theta\in\bC$  is fixed with $\Im\theta >0$.
$1+\chi_2 R_0(z, \theta)W$ and $\chi_2 R_0(z, \theta)\chi_1$ belong to Gevrey class $G^\beta (\Omega_\delta)$ where $\Omega_\delta$ is defined in Proposition \ref{prop6.2}. 
\\

Let $\{\psi_j, j=1, \cdots, m\}, Q, Q'$ be defined as in the proof of Theorem \ref{th6.4}.
 Then $-1$ is not an eigenvalue of compact operator $Q'( G_0W)Q'$.
Since 
 $ Q'(\chi_2 G_0W)Q'$ converges to $Q'G_0WQ'$ in operator norm as $R\to\infty$, $-1$ is not an eigenvalue of $Q'\chi_1 G_0WQ'$
  if $R \ge R_1$ for some $R_1\ge R_0 $, $R_1$  sufficiently large. 
Then $Q'(1+\chi_2 G_0W)Q'$ is invertible on $\range Q'$, so is  $Q'(1+  \chi_2 R_0(z,\theta) W)Q'$ for 
$z \in \Omega_\delta(c) =\Omega_\delta\cap \{z; |z|<c\}$ for some $c>0$.  The inverse 
\be \label{6.31}
E_0(z,\theta) =  (Q'(1+  \chi_2 R_0(z, \theta) W)Q')^{-1}Q'
\ee
is uniformly bounded in $z$  (see Proposition \ref{prop5.5}) and by (\ref{6.36}) it belongs to  Gevrey class $G^\beta(\Omega_\delta(c))$. \\

 Define   $S_1 : \bC^m \to L^2$ and $T_1 : L^2 \to \bC^m$ by 
 \be
 S_1 = \chi_1 S, \quad T_1 = T\chi_1 
 \ee
where $S, T$ are defined in Theorem \ref{th6.4}. By Lemma \ref{lem6.3},
\be
S_1T_1 = Q' + O(e^{-cR^{1-\mu}}), \quad  T_1S_1 = 1  +  O(e^{-cR^{1-\mu}})
\ee
for some $c>0$. Let $W(z,\theta) =1 +\chi_2 R_0(z, \theta) W$. Consider  the Grushin problem 
\be
\left( \begin{array}{cc}
W(z,\theta) & S_1\\
 T_1 &  0
\end{array}\right)  : L^2\otimes \bC^m \to L^2\otimes \bC^m. 
\ee
One has 
\be
\left( \begin{array}{cc}
W(z,\theta) & S_1\\
 T_1 &  0
\end{array}\right) \left( \begin{array}{cc}
E_0(z,\theta) & S_1\\
 T_1 &  - T_1W(z,\theta)S_1
\end{array}\right) = 1 + \vR(z, \theta)
\ee
where
\be
\vR(z, \theta) =  \left( \begin{array}{cc}
QW(z, \theta)E_0(z,\theta) + S_1T_1 - Q & (1-T_1)W(z, \theta)S_1\\
 T_1 E_0(z, \theta) &  T_1 S_1 -1 
\end{array}\right).
\ee
$\vR(z, \theta)$ is  sum of a nilpotent matrix and a matrix of order $O(e^{-cR^{1-\mu}})$.
Hence $1 +\vR(z)$ is invertible  $z\in \Omega_\delta(c)$ if $R>R_1$ is sufficiently large. This proves the Grushin problem is invertible from the right.
Similarly one can show it is invertible from the left, therefore it is invertible with inverse given by 
\be
\left( \begin{array}{cc}
E_0(z,\theta) & S_1\\
 T_1 &  - T_1W(z,\theta)S_1
\end{array}\right) (1+\vR(z))^{-1}  :=  \left( \begin{array}{cc}
E(z) & E_+(z)\\
E_{-}(z)  &  E_{-+}(z) 
\end{array}\right)
\ee
As usual, one has the formula
\begin{equation}
(1+ \chi_2 R_0(z,\theta) W)^{-1} = E(z) - E_+(z) E_{-+}(z)^{-1} E_-(z).
\end{equation}
$E_{-+}(z) $ is of the form
\[
E_{-+}(z)  =  - T_1W(z,\theta)S_1 (1 + O(e^{-cR^{1-\mu}})) + O(|z|^2)
\]
By the choice of $\chi_1,  \chi_2$, one has 
\bea
 T_1W(z,\theta)S_1 &=& T_1 ( 1 + R_0(z,  \theta)W)S_1 \\
 &=& z T_1 G_1(\theta)W)S_1  + O(|z|^2) = z T_1 G_1 WS_1 + O(|z|^2).
 \eea
 By the calculation made in the proof of Theorem \ref{th6.4},  one sees $\Psi_1 = T_1G_1S_1$ is an invertible matrix (if $R$ is large enough).
 Consequently $E_{-+}(z)$ is invertible for $z\in \Omega_\delta(c)$ with inverse of the form. 
 \be
 E_{-+}(z)^{-1}  = -\f{1}{z} \Psi_1 (1 + O(e^{-cR^{1-\mu}})) + B(z)
 \ee
 where $B(z)$ belongs to $G^\sigma(\Omega_\delta(c))$. This proves the existence of an asymptotic expansion for $\chi_1R(z) \chi_1$ for  $z \in \Omega_\delta(c)$ of the form
 \be
 \chi_1 R(z)\chi_1 = -\f{U}{z} + R_2(z)
 \ee
 with $R_2(z)$ satisfying  Gevrey estimates of order $\sigma$ on $\Omega_\delta(c)$. To determine $U$, we remark that
 since $\vA \subset \vV$,   Theorem \ref{th6.4} applied to $R(z)$ with $z\in \Omega_\delta(c)\cap\{\Re z \le 0\}$ gives $U=\chi_1 \Pi_0\chi$.
\ef

{\bf \noindent Proof of Theorem \ref{th1.3}.} Theorem \ref{th1.3} (a) and (b) are respectively deduced from Theorems \ref{th6.4} and \ref{th6.4b} and  the formulas for $t>0$
\bea \label{6.40}
e^{-tH} -\sum_{\lambda \in \sigma_d(H), \Re \la \le 0} e^{-tH}\Pi_\lambda  &=& \f i{2\pi} \lim_{\ep \to 0_+} \int_{\Gamma(\ep)} e^{-tz} R(z) dz + O(e^{-ct}) \\
\label{6.40b}
\chi(e^{-itH} -\sum_{\lambda \in \sigma_d(H)\cap \bR_-} e^{-i tH}\Pi_\lambda)\chi   &=& \f i{2\pi} \lim_{\ep \to 0_+} \int_{\Gamma_\eta(\ep)} e^{-itz}\chi R(z) \chi dz + O(e^{-ct})
\eea
where $c>0$ and
\beas 
\Gamma(\ep) &= & \{z; |z| \ge \ep, \Re z \ge 0, |\Im z| = C (\Re z)^{\mu'}\}\cup\{z; |z|=\ep,  | \arg z| \ge \omega_0\} \\
\Gamma_\eta(\ep) &= & \{z =r e^{-i\eta}, r\ge \ep\}\cup\{z = -re^{i\theta}, r\ge \ep\} \cup\{z;  |z| = \ep,  -\eta \le \arg z \le \pi + \eta\}
\eeas
for some appropriate constants $C, \mu'>0, \eta>0$. In particular, $\eta>0$ is chosen such that $H$ has no eigenvalues with negative imaginary part above $\Gamma_\eta(\ep)$. Here $\omega_0$ is the argument of the point $z_0$ with $|z_0| =\ep$, $\Re z_0>0$ 
 and $\Im z_0 = C (\Re z_0)^{\mu'}$. Remark that the subexponential time-decay estimates are derived from Gevrey estimates of $R_1(z)$ and $R_2(z)$ at zero and  their Taylor expansion of order $N$ with $N$ chosen appropriately in terms of $t>0$. See the proof of Theorem \ref{th4.3} for $e^{-tH_0}$.
\ef

\begin{remark}
As an example of applications of  Theorem \ref{th1.3}, consider the Witten-Laplacian
\be
-\Delta_U = ^t\nabla_U\cdot\nabla_U
\ee
where $\nabla_U =e^{-U}\nabla e^{U}$ and $U\in C^2(\bR^n)$. Then
\[
-\Delta_U = - \Delta + (\nabla U)(x)\cdot (\nabla U)(x) - \Delta U(x)
\]
If  $U \in C^2(\bR^n; \bR)$   satisfies for some $\rho \in ]0, 1[$ and $c_1, C_1 >0$,
 \be \label{assU}
 U(x) \ge c_1\w{x}^\rho, \quad |\nabla U(x)| \ge c_1\w{x}^{\rho -1}, \quad  |\p_x^\alpha U(x)| \le C_1 \w{x}^{\rho -|\alpha|}
\ee 
for $x$ outside some compact and for $\alpha \in \bN^n$ with $|\alpha|\le 2$.  Then  $-\Delta_U$ can be decomposed as $  -\Delta_U = H_0 + W(x) $ where $H_0$ satisfies the conditions of Theorem \ref{th1.1} with $\mu =1-\rho$  and $W(x)$ is of compact support. 
Zero is a simple eigenvalue of  $-\Delta_U$ embedded in its continuous spectrum $[0, +\infty[$.
As consequence of Theorem \ref{th1.3}, one obtains the following result.
Let $\varphi_0 (x)$ be a normalized eigenfunction of $-\Delta_U $ with eigenvalue zero:
\be
\vp_0(x) =Ce^{-U(x)}, \quad \|\vp_0\| =1.
\ee 
Then for any $a>0$, there exist some constants $C_a, c_a >0$
 \be \label{1.41} 
 \| e^{t \Delta_U}f - \w{\varphi_0, f}\varphi_0 \| \le C_a e^{-c_a t^{\f{\rho}{2-\rho}}} \|e^{a\w{x}^\rho}f\|
  \ee
 for $t>0$ and $f$ such that $e^{a\w{x}^\rho}f \in L^2$. 
 Note that the subexponential convergence estimate (\ref{1.41}) without the explicit remainder estimate on $f$ is proved  in \cite{dfg} by method of Markov processes. 
\end{remark}

 \sect{Threshold spectral analysis in non-selfadjoint case}
 
\subsection{The general case}
Finally we study the case  zero is an embedded eigenvalue of the non-selfadjoint Schr\"odinger operator $H$. 
 Let $V_0\in \vV$. Then zero is an eigenvalue of
$H$ if and only if $-1$ is an eigenvalue of compact operator $K=G_0W$ on $L^2(\bR^n)$.  The algebraic multiplicity $m$ of eigenvalue $-1$ of $K$  is finite, although we do not know how to define the algebraic multiplicity of zero eigenvalue of $H$.
Let $\pi_1 : L^{2} \to L^{2}$ be the associated Riesz projection of $K$ defined by :
\[
\pi_1 =\f 1{2\pi i} \int_{ |z+1|=\ep} (z-K)^{-1} dz
\]
for $\ep>0$ small enough. 
 Then
\be
m =\ran \pi_1.
\ee
$\pi_1 $ is continuous on $L^{2, s}$ for any $s\in\bR$ and 
$\pi_1^* : L^{2} \to L^{2}$ is the Riesz projection of $ K^*$ associated with the eigenvalue $-1$.
 \\

By Corollary \ref{cor4.3}, $R_0(z)W$ is continuous in $z \in \Omega(\delta)$, where 
\[
\Omega(\delta) =\{ |z|<\delta; Re z <0 \mbox{ or } \Re z \ge 0 \mbox{ and} |\Im z| > M(\Re z)^{\mu'}\}
\]
 for some $M, \delta$ and $\mu'>0$. Denote $\pi_1' = 1 -\pi_1$. $\pi'(1+G_0W)\pi'$ is injective on the range of $\pi_1'$. The Fredholm Theorem implies that $(\pi'_1 (1+ G_0W) \pi'_1)^{-1}$ is invertible on $L^2$. It follows that
\be
B_1(z) = (\pi'_1 (1+ R_0(z) W) \pi_1')^{-1}\pi_1'
\ee
 is well defined on $\Omega(\delta)$ if $\delta>0$ is sufficiently small. In addition $B_1(z)$ is uniformly bounded there.
 Since $R_0(z)  W$ satisfies Gevrey estimates of order $\sigma$ for $z$ near $0$ with $\Re z <0$ and $|\Im z| < -C \Re z$, $C>0$, 
 so does $B_1(z)$. $\pi_1 (1+ R_0(z) W) \pi_1$ is of finite  rank. Set
 \be \label{ome}
 \omega(z) = \det (\pi_1 (1+ R_0(z) W) \pi_1).
 \ee
Then $\pi_1 (1+ R_0(z) W) \pi_1$ is invertible if and only if $\omega(z)\neq 0$. 
  $\omega(z)$ satisfies the Gevrey estimates of order $\sigma$ at point $z=0$ and has an asymptotic expansion of the form
  \be\label{ome1}
  \omega(z) = \sum_{j=1}^N \omega_j  z^j + O(|z|^{N+1}), z\in \Omega(\delta), 
  \ee
for any $N$. \\

\begin{prop} \label{prop7.1} Assume that 
\be \label{7.0}
\omega(z) =  \omega_k z^k + O(|z|^{k+1})
\ee
 for some $k\in \bN^*$ and $\sigma_k\neq 0$.
Let $-1$ be an eigenvalue of $G_0W$ with algebraic multiplicity $m$.
Then there exist operators  $C_{-j}$, $j=1, \cdots, k$ with rank less than or equal to $m$ such that 
\be \label{7.1}
R(z) = \f{C_{-k}}{z^k}+\cdots +  \f{C_{-1}}{z}  + R_3(z)
\ee
 for $z\in \Omega(\delta)$. The remainder $R_3(z)$ satisfies the estimates:
$\exists C,\mu', \delta >0$ such that 
\be \label{7.2}
\|\w{x}^{-s} R_3(z)\|  + \|R_3(z) \w{x}^{-s} \| \le C_s
\ee
for $s> 2\mu + \f 1 {\mu'}$ and $z\in \Omega(\delta)$.
\end{prop}
\pf Since $\omega_k \neq 0$, $\pi_j (1+ R_0(z)W)\pi_1$ is invertible on the range of $\pi_1$ for $z\in \Omega(\delta)$ with $\delta>0$ small enough. Set $B_0(z) =(\pi_j (1+ R_0(z)W)\pi_1)^{-1}\pi_1$. Then $\omega(z) B_0(z)$ has same continuity properties as $\pi_j (1+ R_0(z)W)\pi_1$
and
\be
B_0(z) = z^{-k}B^{(0)}_{-k} +  \cdots  z^{-1}B^{(0)}_{-1} + R^{(0)}(z)
\ee
where $B^{(0)}_{-j}$, $j=1, \cdots, k$, are operators of rank $\le m$ and $ R^{(0)}(z)$ is uniformly bounded for $z\in\Omega(\delta)$.
One can check that
\be
(1+R_0(z)W)(B_0(z) + B_1(z)) =1 + O(|z|),  (B_0(z) + B_1(z)) (1+R_0(z)W) =1 + O(|z|)
\ee
in $\vL(L^2)$. Therefore $(1+R_0(z)W)$ is invertible for $z\in \Omega(\delta)$ if $\delta>0$ is small enough and
\be
((1+R_0(z)W)^{-1} =  B_0(z) ( 1 + O(|z|))^{-1} + O(1), \quad z \in \Omega(\delta).
\ee
(\ref{7.1}) can be now derived from the equation $R(z) = (1+ R_0(z)W)^{-1} R_0(z)$.
\ef

\begin{remark}\label{rmk7.1}  Under the conditions of Theorem \ref{th6.1.1}, let $\la \in r_+(H)$ be an outgoing positive resonance of $H$.
Let $\pi_1$ denote the Riesz projection of eigenvalue $-1$ of $R_0(\la+i0)W$ as operator on $L^{2,-s}$, $\f 1 2 < s <\rho-\f 1 2$. Denote 
\[
\omega (z) = \det \left(\pi_1 (1+ \chi R_0(z)W)\pi_1\right)
\]
for $z\in \bC_+$, where $\chi \in C_0^\infty$ with $\chi W =W$. Then one can show that $\omega(\la)=0$ and  $\omega(z)\neq 0$ for $z$ with $\Im z >>1$. 
In addition, $\omega(z)$ extends to a holomorphic function into a complex neighbourhood of $\la$. One concludes that
there exist some $k \in \bN^*$ and some $\omega_k \neq 0$ such that 
\be
\omega(z) =  \omega_k (z-\la)^k + O(|z-\la|^{k+1})
\ee
for $z$ in a complex neighbourhood of $\lambda$. This means a condition analogous to  (\ref{7.0}) is satisfied for positive resonances under some analyticity condition on potentials. The proof of Proposition \ref{prop7.1}  and formula (\ref{6.2.4}) allow to conclude that  the meromorphic extension from $\bC_+$  of the cut-off resolvent  $\chi R(z) \chi $  admits an expansion around $\la$ of the form:
\be
\chi R(z) \chi  = \chi (\f{C_{-k}}{(z-\la)^k}+\cdots +  \f{C_{-1}}{z-\la})  + \tilde R_3(z))\chi
\ee
for $z$ near $\la$, where $C_{-j}$ is of rank less than or equal to $m_+(\lambda)$ and $\tilde R_3(z)$ is holomorphic in a neighbourhood of $\lambda$.
Here $m_+(\lambda) =\ran \pi_1$ is the algebraic multiplicity of eigenvalue $-1$ of $R_0(\la+i0)W$.
\end{remark}

\subsection{Representation of the Riesz projection} 
In order to give some more precisions on the resolvent expansion given in Proposition \ref{7.1}, we study in more details the Riesz projection $\pi_1$ accociated with eigenvalue $-1$. Assume from now on that this eigenvalue is geometrically simple. Set $K=G_0W$.
Then
\be
\dim \ker (1+K) =1, \quad \ran \pi_1 =m.
\ee
 Operator $1+K$ being nilpotent on $\range \pi_1$, there exists some function $\phi_m \in \mbox{ range } \pi_1$ such that
\be
\phi_j = (1+K)^{m-j}\phi_m \neq 0, \quad j=1, \cdots, m.
\ee
One has 
\be
(1+K)\phi_1 = 0,  \quad  (1+K)\phi_j = \phi_{j-1}, \quad 2 \le j \le m.
\ee
$\phi_1, \cdots, \phi_m$ are linearly independent. Denote $J$ the operation of complex conjugaison $ J : f \to \overline{f}$.
 Remark that $ H_0^* = JH_0J$, $H^* = J H J$. One has  
 \be
 J W K = K^* \overline {W} J.
 \ee
 It follows that 
 \be 
 J W \pi_1 = \pi_1^*J W. 
 \ee
 Denote
\be
\phi_j^* = \overline{W} \overline{\phi_j}.
\ee
Then
\be
(1+K^*)\phi_1^* = 0,  \quad  (1+K^*)\phi^*_j = \phi^*_{j-1}, \quad 2 \le j \le m.
\ee
Since $\phi_1^* \neq 0$, it follows that $\phi_j^* \neq 0$ for all $1 \le j \le m$.
From this, we deduce that $\{\phi_j^*, j=1,\cdots, m\}$ is linearly independent and that $\ran \pi_1 \le \ran \pi_1^*$. Similarly, 
using the relation
\be
JG_0\pi_1^* = \pi_1 JG_0
\ee
one can prove that $\ran \pi_1 \ge \ran \pi_1^*$, which gives

\begin{lemma} \label{lem6.6} One has
\be
\ran \pi_1 = \ran \pi_1^* = m
\ee
and $JW$ is a bijection from $\range \pi_1$ onto $\range \pi^*_1$.
\end{lemma}

\begin{lemma} \label{lem6.7} The bilinear form $B(\cdot, \cdot)$ defined on $\range \pi_1$ by
\be
B(\varphi, \psi)  =\w{\varphi, J W\psi} = \int_{\bR^n} W(x) \varphi(x)\psi(x) \; dx
\ee
is non-degenerate.
\end{lemma}
\pf Let $\phi \in\range \pi_1$ such that
\[
\int_{\bR^n} W(x) \phi(x) \varphi(x) \; dx =0
\] 
for all $\varphi \in \range \pi_1$. This means that $\phi \in (\range \pi_1^*)^\perp = \ker \pi_1$, which implies that
$\phi = \pi_1 \phi =0$. So $B(\cdot, \cdot)$ is non-degenerate.
\ef

As a consequence of Lemma \ref{lem6.6}, if $m=1$, then eigenfunction $\vp$ of $H$ associated with zero eigenvalue satisfies
\be
\int_{\bR^n} W(x) (\vp(x))^2 \; dx \neq 0.
\ee

\begin{lemma} \label{lem6.8} There exist $\chi_j \in \ker(1+K)^{m-j+1}$, $j=1,\cdots, m$,  such that 
\be
\w{\phi_i, \chi_j^*} = B(\phi_i, \chi_j) =\delta_{ij}, \quad 1\le i, j \le m,
\ee
where  $\chi_j^* = JW \chi_j$, $\delta_{ij} =1$ if $i=j$, and $\delta_{ij} =0$ if $i \neq j$.
\end{lemma}
\pf  We use an induction to prove that for any $1\le l \le m$,  there exist $\varphi_j \in \ker (1+K)^j$, $1\le j \le l$  such that
\be \label{r1}
B(\varphi_i, \phi_{m-j+1}) =\delta_{ij}, \quad  1\le j \le i \le l.
\ee
Since $\phi_1 \in \ker (1+K)$ and $\phi_j^* \in \range (1+K^*)$ for $1\le j \le m-1$, one has $\w{\phi_1, \phi_j^* }=0$ for $j=1, \cdots, m-1$. By lemma \ref{lem6.7}, one has necessarily $c_1= \w{\phi_1, \phi_m^*} \neq 0.$
Set 
\be
\varphi_1 =\f 1{c_1} \phi_1.
\ee
Then $\varphi_1 \in \ker(1+K)$ and  $B(\varphi_1, \phi_m)=1$. (\ref{r1}) is true for $l =1$. Assume now that (\ref{r1}) is true for some $l=k-1$, $2\le k \le m$. Set
\be
\phi'_k =\phi_k -\sum_{j=1}^{k-1} B(\phi_k, \phi_{m-j+1})\varphi_j
\ee
Then $\phi'_k \neq 0$, $\phi'_k \in \ker(1+K)^k$ and
\[
B(\phi'_k, \phi_{m-j+1}) =0, \quad j=1, \cdots, k-1.
\]
Since $\phi'_k \in \ker(1+K)^k$,   one has also 
\be \label{e3}
\w{\phi'_k, \phi_j^*} = B(\phi'_k, \phi_j) =0
\ee
for $j= 1, \cdots,  m-k$, because $\phi_j^* = (1+ K^*)^{m-j}\phi_m^* $ belongs to the range of $(1+K^*)^k$ if $1 \le j \le m-k$.
By Lemma \ref{lem6.7}, the constant $c_k = B(\phi'_k, \phi_{m-k+1})$ must be nonzero. Set
\be
\varphi_k =  \f{1}{c_k}\phi'_{k}.
\ee
Then (\ref{r1})  is proved for $l =k$. By an induction, one can construct 
$\varphi_j$, $1\le j \le m$, such that (\ref{r1}) holds with $l=m$. By (\ref{e3}), one has also
$B(\varphi_i, \phi_{m-j+1})=0$ if $i>j$.  
Lemma \ref{lem6.8} is proved by taking $\chi_k = \varphi_{m-k+1}$, $1 \le k\le m$. \ef

One has the following representation of the Riesz projection $\pi_1$. \\

\begin{cor} \label{cor6.9}
One has 
\be
\pi_1 u = \sum_{j=1}^m \w{u, \chi_j^*} \phi_j, \quad u\in H^{1,-s},  s>1.
\ee
\end{cor}
\pf Denote $\pi$ the operator $\pi : u \to  \sum_{j=1}^m \w{u, \chi_j^*} \phi_j$. Then it is clear that $\pi^2=\pi$ and $\range \pi = \range \pi_1$. 
It is trivial that $\ker \pi_1 \subset \ker \pi$. If $u \in \ker \pi$, then $\w{u, \chi_j^*}=0$ for $j=1, \cdots, m$. Therefore $u \in (\range \pi_1^*)^\perp = \ker \pi_1$ which implies that $\ker \pi \subset \ker \pi_1$. 
This shows that $\ker \pi_1 = \ker \pi$.  This proves $\pi =\pi_1$.
\ef

From the proof of Lemma \ref{lem6.8}, one sees that if $-1$ is a simple eigenvalue of $K$ ($m=1$), then the associated  Riesz projection is given by
\be
\pi_1= \w{\cdot, \vp^*}\varphi
\ee
where $\varphi$ is an eigenfunction of $K$ with eigenvalue $-1$  normalized by
\[
\int_{\bR^n} W(x)(\vp(x))^2 \; dx =1.
\]

\subsection{Resolvent expansion at threshold.}
To study the singularity of the resolvent $R(z)$ at threshold zero, we use the resolvent equation
\[
R(z) = (1+ R_0(z)W)^{-1}R_0(z)
\]
 for $z\not\in \sigma(H)$ and  study the following Grushin problem in $L^2 \times \bC^m$:
\be
\left ( \begin{array}{cc}
1+ R_0(z)W &  S
\\
T & 0
\end{array} \right) : L^{2} \times \bC^m \to L^{2} \times \bC^m
\ee
where 
\bea
S : & &\bC^m \to L^{2}, c=(c_1, \cdots, c_m) \to Sc = \sum_{j=1}^m c_j \phi_j, \\
T : & & L^{2} \to \bC^m, f \to  Tf = (\w{f, \chi^*_1} , \cdots, \w{f, \chi_m^*}). \\
\eea
Then $ST =\pi_1$ and $TS =I_n$. Since $K$ commutes with its Riesz projection $\pi_1$ and since $1+K$ is injective on $\range \pi_1'$ where $\pi_1'=1-\pi_1$, $ 1+ K$  is invertible on the range of $\pi_1'$. By an argument of perturbation, $ \pi_1'(1+ R_0(z)W)\pi_1'$ is invertible
on range of $\pi_1'$ for $z \in \Omega_1(\delta)$ if $\delta>0$ is appropriately small and its inverse $E(z)$ is uniformly bounded
on $\Omega_1(\delta)$ where 
\be
E(z)  = (\pi_1'(1+ R_0(z)W)\pi_1')^{-1}\pi_1'
\ee
By the arguments used in Section 6.1, $E(z)$ belongs to the Gevrey class $G^\sigma(\Omega_1(\delta))$ with $\sigma = 1 +\gamma$. One can check that for $z\in \Omega_1(\delta)$,
\be
\left ( \begin{array}{cc}
1+ R_0(z)W &  S
\\
T & 0
\end{array} \right)^{-1} = \left ( \begin{array}{cc}
E(z) &  E_+(z)
\\
E_-(z) & E_{-+}(z)
\end{array} \right) 
\ee
where 
\bea
E_+(z) &= & (1-  E(z)R_0(z)W)S \\
E_-(z) & = & T(1- R_0(z)WE(z)) \\
E_{-+}(z) & = &  -T(1+ R_0(z)W)S + TR_0(z)WE(z) R_0(z)WS. \label{7.33}
\eea
 It follows that  $z\not\in \sigma(P)$ if and only if 
$\det E_{-+}(z) \neq 0$ and one has
\be \label{Rz}
(1+R_0(z)W)^{-1} = E(z) - E_+(z)  E_{-+}(z)^{-1} E_-(z).
\ee
By operations on Gevrey functions ((\ref{6.34}) -(\ref{6.37})), $E_{-+}(z)$ is $m\times m$-matrix valued Gevrey function for $z \in \Omega_1(\delta)$.  Therefore it has an asymptotic expansion near $0$ up to any order
\be
E_{-+}(z) = B_0 + B_1 z + \cdots + B_N z^n + O(|z|^{N+1})
\ee 
where $B_j$ is some $m\times m$ matrix. More precisely, since $ TR_0(z)WE(z) R_0(z)WS = O(|z|^2)$, $E_{-+}(z)$ verifies
\bea
E_{-+}(z) & = & \left(-\w{ (1+R_0(z)W)\phi_k, \chi_j^*} \right)_{1\le j, k \le m } + O(|z|^2) \\
 & = & -\left( \begin{array}{cccc}  
 0&1 & \cdots & 0 \\
 0&0 & \ddots & \vdots \\
 \vdots &\vdots & \ddots  &1 \\
 0 & 0 & \cdots  & 0
 \end{array}  \right)   - z  \left( \begin{array}{cccc }  
  b_{11}&\cdots & \cdots & b_{1m}\\
  \vdots& \ddots &   & \vdots \\
 \vdots &   &\ddots  & \vdots \\
 b_{m1}  & \cdots &\cdots & b_{mm}
 \end{array}  \right) + O(|z|^2) \nonumber
\eea
where 
\be
b_{jk} = \w{G_1 W\phi_k, \chi_j^*}.
\ee
Note that $\phi_1$ and $\chi_m$ belong to $\ker(1+G_0W)$  and $\chi_m^* = J W\chi_m$, they are rapidly decreasing, by Lemma \ref{lem6.3}.
One can calculate
\beas
b_{m1} &=&  \lim_{\la \to 0_- } \w{\f 1 \la (1+ R_0(\la)W) \phi_1,  J W\chi_m} \\
& = &  -\lim_{\la\to 0_-} \w{R_0(\la) \phi_1,  J W\chi_m}   \\
& = &   -\lim_{\la\to 0_-} \w{ \phi_1,  J R_0(\la) W\chi_m}\\
&= & \w{\phi_1, J \chi_m}.
\eeas
Similarly one can calculate for $2 \le j \le m$
\[
b_{mj} =  - \w{ W\phi_j,  J G_0\chi_m} =\w{\phi_{j}-\phi_{j-1}, J \chi_m}.
\]
Summing up, we have proved the following
\\

\begin{prop} \label{prop7.6} $\det E_{-+}(z)$ is a Gevrey function of order $\sigma$ for $z \in \Omega_1(\delta)$ 
and has an asymptotic expansion in powers of $z$ 
\be
\det E_{-+}(z) = \sigma_1 z + \cdots \sigma_N z^N + O(|z|^{N+1})
\ee
for any $N$, where 
\be
\sigma_1 = - b_{m1}.
\ee
\end{prop}

\begin{theorem}\label{th6.5} Let $H_0 = -\Delta + V_0(x)$ and $H = H_0 + W(x)$
with $V_0\in \vV$ and $W\in L^\infty_{\rm comp}$.  Assume that zero eigenvalue of $H$
is geometrically simple.
\\

(a).  Suppose that
\be\label{7.61}
\det E_{-+}(z) = \sigma_k z^k +  O(|z|^{k+1})
\ee
for some $\sigma_k \neq 0$, $k \ge 1$. Then there exist operators  $C_j$, $j=-k, \cdots, -1$ with ranks less than or equal to $m$ such that 
\be
R(z) = \f{C_{-k}}{z^k}+\cdots +  \f{C_{-1}}{z}  + R_3(z)
\ee
 for $z\in \Omega_1(\delta)$, where $C_{-j}$, $ 1 \le j \le k-1$, are of rank less than or equal to $m$ and $C_{-k}$ is a rank one operator given by
 \be
 C_{-k} =\w{\cdot, J\vp_1}\vp_1,
 \ee
with $\vp_1$ an eigenfunction of $H$ associated with zero eigenvalue.
 The remainder $R_3(z)$ satisfies the estimates:
$\exists C, \mu', \delta >0$ such that 
\be \label{6.52}
\|\w{x}^{-s} R_3(z)\|  + \|R_3(z) \w{x}^{-s} \| \le C_s
\ee
for $s> 2\mu + \f 1 {\mu'}$ and $z\in \Omega_1(\delta)$; and  for any  $ a>0$,  $\exists C_a, c_a >0$ such that
\be \label{6.53}
\| e^{-a\w{x}^{1-\mu}} R_3^{(N)}(z)\| + \|   R_3^{(N)}(z)e^{-a\w{x}^{1-\mu} }\| \le C_a c_a^N N^{\sigma N},
\ee
for any $N \in \bN$ and $z \in \Omega_-$ where $\Omega_1(\delta)$  and $\Omega_-$ are the same as in Theorem \ref{th6.4}.
\\

(b). Suppose in addition that there exists an eigenfunction $\vp_0$ of $H$ associated with eigenvalue zero such that 
\be \label{6.50}
\int_{\bR^n} (\vp_0(x))^2 dx =1.
\ee 
Then Condition (\ref{7.61}) is satisfied with $k=1$ and one has
\be
C_{-1} = -\w{\cdot, J\vp_0}\vp_0.
\ee
\end{theorem}
\pf (a). The existence of the resolvent expansion is proved in Proposition \ref{prop7.1} and the Gevrey estimates of the remainder can be obtained in the same way as in Theorem \ref{th6.4}. We only calculate $C_{-k}$. 
Under the condition (\ref{7.61}), one has
\be
E_{-+}(z)^{-1} = \f{ ^t {\rm Com } \; E_{-+}(z) }{\det E_{-+}(z) } = z^{-k} C  + O(|z|^{-k+1}) 
\ee
for $z\in \Omega_1(\delta)$, where 
\[
C= \left( \begin{array}{cccc}  
 0&    \cdots  & 0& \sigma_k^{-1}\\
 \vdots & \ddots &    &  0 \\
 \vdots & &  \ddots & \vdots \\
 0 &\cdots & \cdots & 0\\
 \end{array}  \right).
 \]
From (\ref{Rz}), one obtains
\be
(1+ R_0(z)W)^{-1} = -z^{-k} SC T  + O(|z|^{-k+1}).
\ee
Using the definition of $S$ and $T$, one sees 
\be
SCT f  = \f 1{\sigma_k} \w{f, \chi_m^*})\phi_1
\ee
Noticing that 
\[
\w{G_0f, \chi^*_m} = \w{f,  G_0^*JW\chi_m} = \w{f, J  G_0W\chi_m} = -\w{f, J\chi_m}, 
\]
 we deduce  from (\ref{R}) that
\be \label{7.40}
R(z) = \f{C_{-k}}{z^k} + O(|z|^{-k+1})
\ee
for $z\in \Omega_1(\delta)$, where $C_{-k} $ is of rank one,  given by
\be
C_{-k} f = \f 1{\sigma_k} \w{f, J\chi_m}\phi_1
\ee
Since $\chi_m$ and $\phi_1$ belong to  the one dimensional space $\ker(1+K)$, $C_{-k}$ can written as
\be
C_{-k}f = \w{f, J\vp_1})\vp_1
\ee
where $\vp_1$ is an eigenfunction of $H$ with eigenvalue zero. This proves part (a).
\\

If (\ref{6.50}) is satisfied, then one has $\chi_m = d_1\vp_0$ and $\phi_1 =d_2 \vp_0$ for some constants $d_j \neq 0$. Therefore
\[
\sigma_1 = -\w{\phi_1, J\chi_m} =- d_1d_2 \neq 0.
\]
Condition (\ref{7.61}) is satisfied with $k = 1$.   Set
\be \label{7.43}
\psi_0 = \sqrt{ \f {d_1}{b_{m1}} } \phi_1.
\ee
Then  $C_ {-1} = -\w{\cdot, J\psi_0}\psi_0$. $\psi_0$ is an eigenfunction of $H$ with eigenvalue zero and
\[
\int_{\bR^n}(\psi_0(x))^2 \; dx = \f{\w{\phi_1, J\chi_m}}{b_{m1}}=1.
\]
Since  zero eigenvalue of $H$ is geometrically simple, one has  $\psi_0= \pm \vp_0$. This proves 
\be \label{7.44}
C_{-1} =  -\w{\cdot, J\vp_0}\vp_0.
\ee 
\ef

\bigskip

\begin{remark} \label{rmk7.2} The methods used here can be applied to other threshold spectral problems. In particular for non-selfadjoint Schr\"odinger operator $H= -\Delta + V(x)$ with a quickly decreasing complex potential $V(x)$ on $\bR^3$:
\be
|V(x)| \le C \w{x}^{-\rho}, \rho>2,
\ee
 our method allows to calculate the low-energy asymptotics of the resolvent $(H-z)^{-1}$ if  zero is a resonance but not an eigenvalue. In fact  using the same reduction scheme and similar  calculations, one can  show in this case $E_{-+}(z)$ takes the form
\bea
\lefteqn{E_{-+}(z) } \nonumber \\
 & = & -\left( \begin{array}{cccc}  
 0&1 & \cdots & 0 \\
 0&0 & \ddots & \vdots \\
 \vdots &\vdots & \ddots  &1 \\
 0 & 0 & \cdots  & 0
 \end{array}  \right)   - z^{\f 1 2}  \left( \begin{array}{cccc }  
  b_{11}&\cdots & \cdots & b_{1m}\\
  \vdots& \ddots &   & \vdots \\
 \vdots &   &\ddots  & \vdots \\
 b_{m1}  & \cdots &\cdots & b_{mm}
 \end{array}  \right) + O(|z|^{\f 1 2 + \ep}). \nonumber
\eea
  The characterization of resonant state  ensures that $b_{m1} \neq 0$. See \cite{jk} in the selfadjoint case. Therefore one can explicitly calculate the leading term of the asymptotic expansion of $(H-z)^{-1}$ for $z$ near $0$ in the case zero is a resonance but not an eigenvalue.
\end{remark}

\bigskip

\begin{theorem}\label{th6.5b} $V_0\in \vA$.   Assume that zero
is a geometrically simple eigenvalue of $H$.  Let $\chi\in C_0^\infty(\bR^n)$ and $\Omega_\delta(c) $ be defined as in Theorem \ref{th6.4b}
Under the condition (\ref{7.61}), the meromorphic extension of $\chi R(z)\chi$ from $\bC_+$ verifies
\be \label{6.51b}
\chi R(z)\chi = \chi (\f{C_{-k}}{z^k}+\cdots +  \f{C_{-1}}{z}  + R_4(z))\chi
\ee
 for $z\in \Omega_\delta(c)$, where $C_{-j}$ is the smae as in Theorem \ref{th6.5} and
 the remainder $R_4(z)$ is continuous up to $z=0$ and satisfies the Gevrey estimates
\be \label{6.6b}
\|\chi R_4^{(N)}(z) \chi\| \le C_\chi C^N N^{\sigma N}
\ee
for $z \in \Omega_\delta(c) \cup\{0\}$. 
In addition if (\ref{6.50}) is true, Condition (\ref{7.61}) is satisfied with  $k=1$ and (\ref{6.51b}) holds with $
C_{-1} = -\w{\cdot, J\vp_0}\vp_0. $
\end{theorem}

Theorem \ref{th6.5b} is derived by combining methods used in  Theorem \ref{th6.4b} and Theorem \ref{th6.5}. The details are omitted. 
\\

\begin{remark} \label{rmk7.3} If zero eigenvalue of $H$ is not geometrically simple, combining  methods used in Proposition \ref{prop7.1} and Theorems \ref{th6.4} and
\ref{th6.4b}, one can show that the resolvent expansions given in Theorems \ref{th6.5} and
\ref{th6.5b} still hold with the corresponding Gevrey estimates on  remainders. But in this case  we can only affirm that $C_{-j}$ is of 
rank $\le m$ for $j=1, \cdots, k$, as in Proposition \ref{prop7.1}.
\end{remark}

{\bf \noindent Proof of  Theorem \ref{th1.4}.} Theorem \ref{th6.5b} implies that outgoing positive resonances are absent in neighbourhood of zero. Therefor under the conditions of this theorem, $r_+(H)$ is at most a finite set. If zero eigenvalue is geometrically simple, the results of Theorem \ref{th1.4} for $e^{-tH}$ can be derived from Theorem \ref{th6.5} and formula (\ref{6.40}) when $V_0\in \vV$ and those for $e^{-itH}$ are obtained from Theorem \ref{th6.5b} and formula (\ref{6.40b} when $V_0\in \vA$. Taking notice of Remark \ref{rmk7.3}, one can prove in the same way the results of
Theorem \ref{th1.4} when zero eigenvalue is not geometrically simple.
\ef

\begin{example}
Consider the non-selfadjoint Witten Laplacian
\[
-\Delta_U = - \Delta + (\nabla U)(x)\cdot (\nabla U)(x) - \Delta U(x)
\]
where  $U\in C^2(\bR^n; \bC)$. Set  $U(x) =U_1(x) + iU_2(x)$ with $U_1, U_2$ real valued functions.
 Assume that $U_1 $ satisfies the condition (\ref{assU}) with $U$ replaced by $U_1$ and that $U_2$ is of compact support with $\| \p_x^\alpha U_2\|_{L^\infty} $ sufficiently small  for $|\alpha | \le 2$. 
Considering $-\Delta_U$ as a perturbation of $-\Delta_{U_1}$,  one can show that $-\Delta_U$ has only one eigenvalue
in a neighbourhood of zero which is in addition geometrically simple. Therefore the eigenfunctions associated zero eigenvalue of $-\Delta_U$ are of the form  $c e^{-U(x)}$ for  some $c\neq 0$, one concludes that the condition (\ref{6.50}) is satisfied if $\| U_2\|_{L^\infty}$ is sufficiently small.
\end{example}


\begin{thebibliography}{10}

\bibitem{ag} S. Agmon,   Spectral properties of Schr\"odinger operators and scattering theory. Ann. Scuola Norm. Sup.
Pisa Cl. Sci. (4) (1975),  no. 2,  151-218.


\bibitem{ac} J. Aguilar, J.M. Combes, A class of analytic perturbations for one-body Schr\"odinger Hamiltonians. Comm. Math. Phys. 22 (1971), 269-279.

\bibitem{bo} D. Boll\'e, Schr\"odinger operators at threshold, 
pp. 173-196, in {\it Ideas and Methods in Quantum and Statistical Physics},
Cambridge Univ. Press, Cambridge, 1992.


\bibitem{pc} P. Cattiaux, Long time behavior of Markov processes, ESAIM Proc. Vol 44, pp. 110-129, 2014.

\bibitem{dfg} R. Douc, G. Fort, A. Guillin,  Subgeometric rates of convergence of $f$-ergodic strong Markov processes,
Stochastic Process.  Appl., 19(2009), 897-923.


\bibitem{FS} S. Fournais, E. Skibsted, Zero energy asymptotics of the resolvent for a
  class of slowly decaying potentials, Math. Z. {\bf 248}
 (2004), 593--633.
 
\bibitem{god} M. Goldberg, A Dispersive Bound for Three-Dimensional Schrödinger Operators with Zero Energy
Eigenvalues, Comm. PDE 35 (2010), 1610-1634.
 
\bibitem{hm} B. Helffer, A. Martinez, Comparaison entre les diverses notions de résonances. (French) [Comparison among the various notions of resonance] Helv. Phys. Acta 60 (1987), no. 8, 992-1003. 





\bibitem{hbst} I. Herbst, Spectral theory of the operator $ (p^2 + m^2)^{\f 1 2 } -\f{Ze^2}{r}$, Commun. in Math. Phys., 53(1977)(3),  285-294.

\bibitem{hun} W. Hunziker, Distorsion analyticity and molecular resonance curves, Annales de
l'I.H.P., (Section Physique Théorique),  Tome 45 (1986) no. 4 , p. 339-358.

\bibitem{jk} A. Jensen  and T. Kato,  Spectral properties of Schr\"odinger operators and time decay of wave 
 functions, Duke Math. J., 46(1979), 583-611.
 
\bibitem{jk} D. Jerison and C.E. Kenig, Unique continuation and absence of positive eigenvalues for
Schr\"odinger operators, Ann. of Math. (2) 121, no. 3 (1985),463-494, With an appendix by
E. M. Stein.

 
 \bibitem{ky} T. Kako, K.  Yajima,  Spectral and scattering theory for a class of non-selfadjoint operators.  Sci. Papers College Gen. Ed. Univ. Tokyo  26  (1976), no. 2, 73-89

 \bibitem{k}T. Kato, Perturbation Theory of Linear Operators, Springer, Berlin, 1980.

\bibitem{kr}  M. Klein, J. Rama, Almost exponential decay of quantum resonance states and Paley-Wiener type estimates in Gevrey spaces,  Ann. Henri Poincaré, 11(2010), 499-537.


\bibitem{lz} T. Li, Z. Zhang, Large time behaviour for the Fokker-Planck equation with general potential, preprint 2016, to appear in Sci. China, Mathematics. 

\bibitem{n} S. Nakamura, Low energy asymptotics for Schr\"odinger operators with slowly decreasing potentials, Commun. in Math. Phys., 161(1994), 63-76.

\bibitem{roy0} J. Royer, Analyse haute-fréquence de l'équation de Helmholtz dissipative, thèse de Doctorat, Univ. Nantes,  décembre 2010. 


\bibitem{roy} J. Royer, Limiting absorption principle for the dissipative Helmholtz equation, Commun. in PDE,  35(8) (2010), 1458-1489.

\bibitem{sai} Y. Saito, The principle of limiting absorption for the nonselfadjoint Schrödinger operator in $R\sp{N}(N\neq 2)$.  Publ. Res. Inst. Math. Sci.  9  (1973/74), 397-428.


\bibitem{sch} J. Schwartz, Some non-selfadjoint operators, Comm. Pure Appl. Math., Vol. XIII (1960), 609-639.

\bibitem{si} B. Simon, Resonances and complex scaling: a rigorous overview. Int. J. Quant. Chem. 14 (1978), 529-542 

\bibitem{sw} E. Skibsted, X. P. Wang,  Two-body threshold spectral analysis, the critical case.  J. Funct. Analysis, 260(6) (2011),   1766-1794.


\bibitem{w1} X.P. Wang,  Asymptotic expansion in time of the Schr\"{o}dinger group on conical manifolds, Ann. Inst. Fourier, Grenoble 56(2006),  1903-1945.

\bibitem{w2} X. P. Wang, Time-decay of semigroups generated by dissipative Schrödinger operators. J. Differential Equations 253 (2012), no. 12, 3523-3542.

\bibitem{w3} X. P. Wang, Large-time asymptotics of solutions to the Kramers-Fokker-Planck equation with a short-range potential. Comm. Math. Phys. 336 (2015), no. 3, 1435-1471.

\bibitem{w4} X. P. Wang, Gevrey type resolvent estimates at the threshold for a class of non-selfadjoint Schr\"odinger operators,  Bruno Pini Mathematical Analysis Seminar, 2015,  69-85.

\bibitem{y1} D. Yafaev,  The low-energy scattering for slowly decreasing potentials, Commun. Math. Phys., 85(1982), 177-196.

\bibitem{y2} D. Yafaev, Spectral properties of the Schr\"odinger operator with a potential having a slow falloff, Funct. Anal. Appli.,  16(1983), 280-286.
\end{thebibliography}
\end{document}